\documentclass[12pt]{article}

\usepackage{amsmath,amsfonts}
\usepackage{fullpage,graphicx}
\usepackage{color}

\newcommand{\B}{\mathbb{B}}
\newcommand{\R}{\mathbb{R}}
\newcommand{\Z}{\mathbb{Z}}

\newcommand{\re}{\mathrm{Re}}

\newcommand{\cP}{\mathcal{P}}
\newcommand{\cX}{\mathcal{X}}
\newtheorem{theorem}{Theorem}[section]
\newtheorem{defn}{Definition}[section]
\newtheorem{asmp}{Assumption}[section]
\newtheorem{corol}[theorem]{Corollary}
\newtheorem{lemma}[theorem]{Lemma}

\newcommand{\stab}{\mathrm{stab}}
\newcommand{\unst}{\mathrm{unst}}
\newcommand{\bif}{\mathrm{bif}}
\newcommand{\fut}{\mathrm{fut}}

\newcommand{\proof}{\noindent {\bf Proof:} }

\newcommand{\qed}{\hfill {\bf QED}\newline}

\begin{document}

\title{Parameter shifts for nonautonomous systems in low dimension: Bifurcation- and Rate-induced tipping}
%\title{
%Parameter shifts for non-autonomous systems: criteria for bifurcation-
%and rate-induced tipping in one dimension
%}

\author{Peter Ashwin\thanks{Centre for Systems, Dynamics and Control, Department of Mathematics, Harrison Building, University of Exeter, Exeter EX4 4QF}, Clare Perryman${}^{*,\dag}$ and Sebastian Wieczorek\thanks{Department of Applied Mathematics, University College Cork,  Western Road, Cork, Ireland}}

\maketitle

\begin{abstract}
  We discuss the nonlinear phenomena of irreversible tipping for
  non-autonomous systems where time-varying inputs correspond to a
  smooth ``parameter shift'' from one asymptotic value to another.  We
  express tipping in terms of properties of local pullback attractors and present some
  results on how nontrivial dynamics for non-autonomous systems can be
  deduced from analysis of the bifurcation diagram for an associated
  autonomous system where parameters are fixed.  In particular, we
  show that there is a unique local pullback point attractor
  associated with each linearly stable equilibrium for the past
  limit. If there is a smooth stable branch of equilibria over the
  range of values of the parameter shift, the pullback attractor will
  remain close to (track) this branch for small enough rates, though
  larger rates may lead to rate-induced tipping. More generally, we
  show that one can track certain stable paths that go along
   several stable branches by pseudo-orbits of the system, for
  small enough rates.  For these local pullback point attractors, we
  define notions of bifurcation-induced and irreversible
  rate-induced tipping of the non-autonomous system.  In
  one-dimension, we introduce the notion of forward basin stability and use this
  to give a number of sufficient conditions for the
  presence or absence of rate-induced tipping. We
  apply our  results to give criteria for irreversible
  rate-induced tipping in a conceptual climate model example.
\end{abstract}

\section{Introduction}

In common language, a {\em tipping point} or {\em critical transition}
is a sudden, large and irreversible change in output of a complex
system in response to a small change of input. There have been a range
of papers in the applied sciences that use the notion of a tipping
point in applications such as climate systems
\cite{Dakos_etal_2008,Lenton_etal_2008,Wieczorek_etal_2010} or
ecosystems \cite{Scheffer2008,Drake_etal_2010}. Further work has
attempted to find predictors or {\em early warning signals} in terms
of changes in noise properties near a tipping point
\cite{Dakos_etal_2008,Ditlevsen_etal_2010,Thompson_Sieber_2010b,Shi_Li_2016,Ritchie2015}, other recent work on tipping includes
\cite{NeneZaikin2012,NishkawaOtt2014}. Indeed there is a large
literature on catastrophe theory (for example \cite{ARN94} and
references therein) that addresses related questions in the setting of
gradient systems.

There seems to be no universally agreed rigorous mathematical
definition of what a tipping point is, other than it is a type of
bifurcation point, though this excludes both
``noise-induced''~\cite{Sutera1981,Ditlevsen_etal_2010} and
``rate-induced''~\cite{Scheffer2008,Wieczorek_etal_2010,PerryWiecz2014}
transitions that have been implicated in this sort of phenomenon.
Other work on critical transitions includes for example~\cite{Kuehn}
that examines the interaction of noise and bifurcation-related tipping
in the framework of stochastically perturbed fast-slow systems. In an
attempt to understand some of these issues, a recent paper
\cite{AshWieVitCox2012} classified tipping points into three distinct
mechanisms: bifurcation-induced (B-tipping), noise-induced (N-tipping)
and rate-induced (R-tipping). As noted in \cite{AshWieVitCox2012}, a tipping point
in a real non-autonomous system will typically be a mixture of these
effects but even in idealized cases it is a challenge to come up with
a mathematically rigorous and testable definition of ``tipping point''.

The aims of this paper are:
\begin{itemize}
\item[(a)] to suggest some definitions for B-tipping and
  R-tipping appropriate for asymptotically constant ``parameter
  shifts'' in terms of pullback attractors for non-autonomous systems.
\item[(b)] to show how these definitions allow us to exploit the
  bifurcation diagram of the associated autonomous (i.e. quasi-static) system to obtain a
  lot of information about B- and R-tipping of the non-autonomous
  system associated with a parameter shift.
\item[(c)] to give some illustrative examples of properties that allow
  us to predict or rule out B- and R-tipping in various cases.
\end{itemize}
Restricting to simple attractors and asymptotically constant ``parameter
shifts'' allows us to obtain a number of rigorous results for a natural class of
non-stationary parameter changes that may be encountered in natural sciences and
engineering.

The paper is organized as follows: In Section~\ref{sec:defs} we
consider arbitrary dimensional non-autonomous systems and introduce a
class of parameter shifts $\cP(\lambda_-,\lambda_+)$ that represent
smooth (but not necessarily monotonic) changes between asymptotic
values $\lambda_{-}$ and $\lambda_{+}$. For such parameter shifts, we
relate stability properties of branches of solutions in the autonomous
system with fixed parameters to properties of local pullback point
attractors of the non-autonomous system with time-varying
parameters. This means that the systems we consider are asymptotically
autonomous \cite{Mischaikow,Rasmussen2008} both in forward and
backward time.

After some discussion of (local) pullback attractors, we show in
Theorem~\ref{thm:pullback} that one can associate a unique pullback
point attractor with each linearly stable equilibrium $X_-$ for the
past limit system (i.e. the autonomous system with
  $\lambda=\lambda_-$).  Then, we define two types of tracking: {\em
  $\epsilon$-close tracking} and {\em end-point tracking}. We show in
Lemma~\ref{lem:pbtrack} that if there is a stable branch of equilibria
with no bifurcation points from this $X_-$ to some stable equilibrium
$X_+$ for the future limit system (i.e. the autonomous system with
$\lambda=\lambda_+$), then this unique pullback point attractor will
$\epsilon$-closely track a path of stable equilibria for sufficiently
slow rates.

On the other hand, there may be a number of stable
paths connecting $X_-$ to $X_+$ - each of these with an image that is a
union of stable branches, but possibly
passing through bifurcation points. In this case, the pullback
point attractor can clearly only track at most one of these
paths. Allowing arbitrarily small perturbations along
  the path, there will be pseudo-orbits that can track any stable
  path. In Theorem~\ref{thm:adiabatic} we give such result by stating
necessary conditions for $\epsilon$-close tracking by pseudo-orbits.

Section~\ref{sec:BRtipping} gives some definitions for 
B- and R-tipping that correspond to lack of end-point tracking for parameter shifts.  In
that section we prove several results for one-dimensional systems,
including the following:
\begin{itemize}
\item Certain bifurcation diagrams or parameter ranges will not permit
  irreversible R-tipping (Theorem~\ref{thm:rtipping}, case 1.).
\item Other bifurcation diagrams will give irreversible
  R-tipping for certain parameter shifts (Theorem~\ref{thm:rtipping},
  cases 2. and 3.).
\item There are restrictions on the regions in parameter space where
  irreversible R-tipping can occur
  (Section~\ref{sec:restrictions}).
\end{itemize}
These results use some properties of the paths of attractors. 
We introduce the notion of {\em forward basin stability} where at each point in time the basin of attraction of an attractor for the associated autonomous system contains all earlier attractors on the path.  We prove that:
\begin{itemize}
\item Forward basin stability guarantees tracking (Theorem~\ref{thm:rtipping}
  case 1.).
\item Lack of forward basin stability, plus some additional assumptions, gives testable criteria for irreversible R-tipping
  (Theorem~\ref{thm:rtipping} cases 2. and 3.).
\end{itemize}
We also give an example application to a class of global energy-balance
climate models in Section~\ref{sec:examples}.

In Section~\ref{sec:discuss} we discuss some of the challenges
to extend the rigorous results to higher-dimensional state or
parameter spaces, to cases where noise is present and to parameter
paths that are more general than parameter shifts.

\subsection{The setting: non-autonomous nonlinear systems}

Suppose we have a non-autonomous system of the form
\begin{equation}
\label{eq:ode}
\frac{dx}{dt}=f(x,\Lambda(r t)),
\end{equation} 
where $x\in\R^n$, $\Lambda(s)\in C^2(\R,\R^d)$, $f\in
C^2(\R^{n+d},\R^n)$, $t\in\R$ and the ``rate'' $r>0$ is fixed. For
$x_0\in\R^n$, $t,s\in\R$ and starting at $x(s)=x_0$ we write 
$$
\Phi_{t,s}(x_0):=x(t)
$$
to denote the {\em solution cocycle} \cite{KloedenRasmussen}, where $x(t)$ is the corresponding solution of (\ref{eq:ode}).
We also assume $\Phi_{t,s}(x_0)$ is defined
for all $t>s$ and $x_0$, and note that it satisfies a cocycle equation
$$
\Phi_{t,s}\circ\Phi_{s,u}(x_0)=\Phi_{t,u}(x_0)
$$
for any $t>s>u$ and $x_0\in\R^{n}$.  Note that, in addition to $x_0$, $t$ and $s$,  $\Phi_{t,s}$
depends on the shape of the parameter shift $\Lambda$, its rate of
change $r$, as well as on the function $f$ that defines the
system. The key difference from the evolution operator of an
autonomous system is that a cocycle
depends on both the time elapsed $(t-s)$ and the initial time $s$.

The model (\ref{eq:ode}) can be extended to include stochastic
effects, e.g. by considering a stochastic differential equation of the
form
\begin{equation}
\label{eq:sde}
dx=f(x,\Lambda(r t))dt+ \eta dW_t
\end{equation}
where $W_t$ represents an $n$-dimensional Wiener process and $\eta$ represents a noise amplitude. A random dynamical
systems approach \cite{ArnoldRDS} gives a solution cocycle
that depends also on the choice of noise path.

One can understand a lot about the asymptotic behaviour of
(\ref{eq:ode}) [or (\ref{eq:sde})] from the bifurcation/attractor diagram of the
associated autonomous system
\begin{equation}
\label{eq:odeaut}
\frac{dx}{dt}=f(x,\lambda),
\end{equation}
where $\lambda$ is fixed rather than time-dependent. General
properties of dynamical systems suggest that if $r$ is ``small
enough'' then we can expect, for an appropriate definition of attractor, 
that trajectories of (\ref{eq:ode}) will closely track a branch of attractors of
(\ref{eq:odeaut}) if the branch continues to be linearly stable over the range of the
varying $\lambda$. In this paper we consider only the case of linearly stable
equilbrium attractors for (\ref{eq:odeaut}) where such statements can be made rigorous through use of geometric
singular perturbation theory \cite{FEN79} or uniform
hyperbolicity \cite{BishnaniMackay2002}. The control theory
literature discusses related questions; see for example
\cite{HinPri2011} that allows $\Lambda$ to depend on state as well.
However, close tracking of a continuous path that traverses
several stable branches of (\ref{eq:ode}) cannot always be guaranteed,
even if $r$ is small enough.

If the variation of $\Lambda$ along a branch of attractors brings the
system to a bifurcation point where no stable branches are nearby, we
have the ingredients of a bifurcation-induced or B-tipping point; no
matter how slowly we change the parameter, there will be a sudden and
irreversible change in the state of the system $x$ on, or near where
the parameter passed through a dynamic bifurcation
\cite{BifDelay,DynBifs,NishkawaOtt2014}.

However, B-tipping is not the only way to get a sudden irreversible
change. Even if there is a branch of attractors available for the
system to track, a rate-induced tipping or R-tipping
\cite{AshWieVitCox2012,Wieczorek_etal_2010} can occur; if there is a
critical value of the rate $r$ beyond which the system cannot track
the branch of attractors then it may suddenly move to a different
state. This happens even though there is no loss of stability in the
autonomous system. Nonetheless, the bifurcation diagram of the
associated autonomous system may give us a lot of information and
indeed constraints on when R-tipping may or may not happen.
Nontrivial behaviour that can appear as a result has been recently
studied in \cite{Wieczorek_etal_2010,Mitry_etal_2013,PerryWiecz2014}
in terms of canard trajectories.

%%%%%%%%%%%%%%%%%%%%%%%%%%%%%%%%%%%%%%%%%%%%%%%%%%

\section{Tracking of attractors for non-autonomous systems}
\label{sec:defs}

\begin{figure}%
\centerline{\includegraphics[width=8cm]{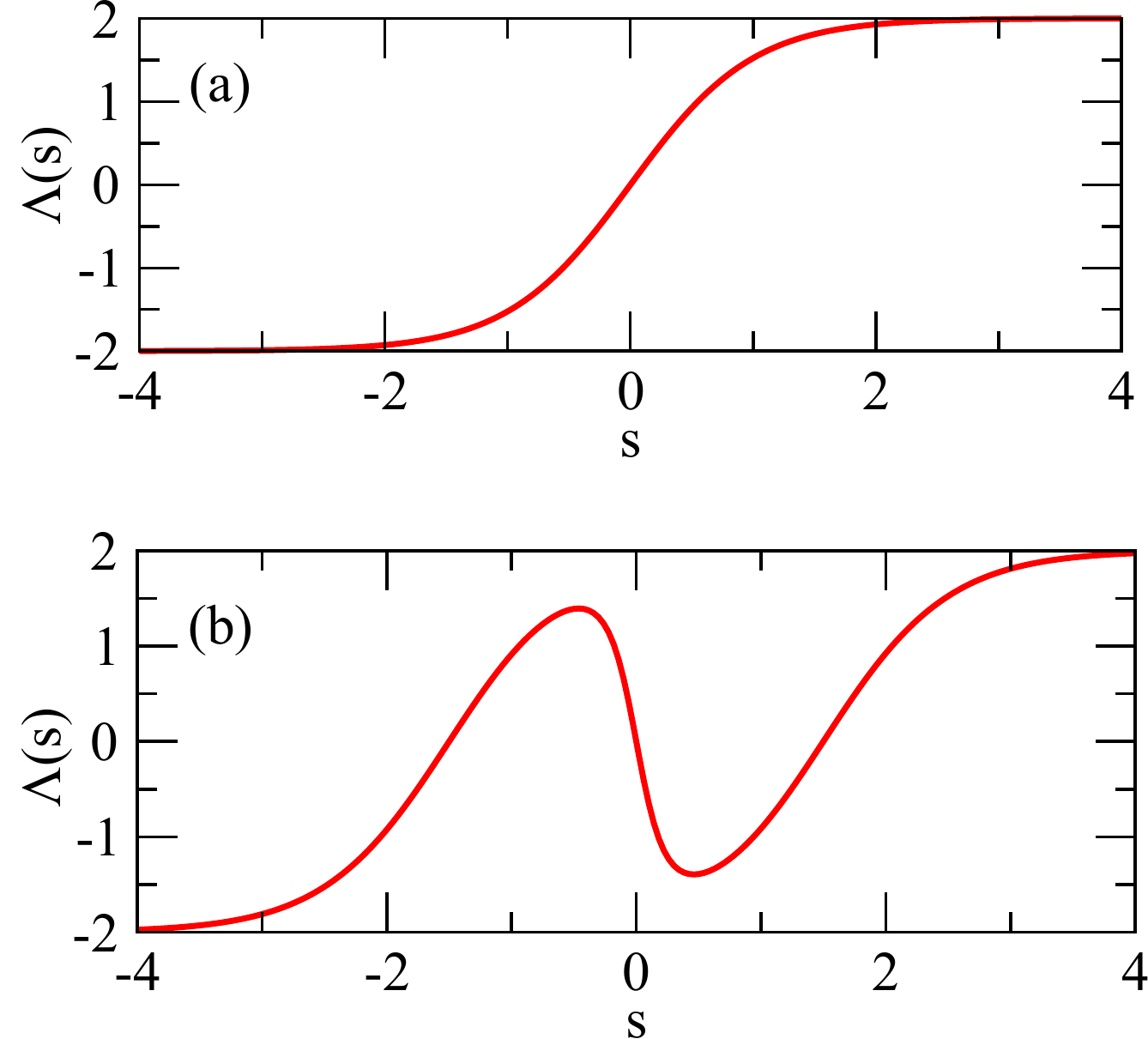}}
\caption{Examples of $\Lambda(s)\in \cP(-2,2)$. (a) shifts monotonically from $-2$ to $2$ while (b) does this shift non-monotonically.
}%
\label{fig:lambdas}
\end{figure}

We consider the behaviour of (\ref{eq:ode}) for one dimensional ($d=1$) 
``parameter
shifts'', i.e. choices of $\Lambda(s)$ that are $C^2$-smooth (but not
necessarily monotonic), bounded and satisfying
$$
\lambda_-<\Lambda(s)<\lambda_+ \mbox{ for each (finite) value of } s,
$$
that limit to the ends of the interval, $\Lambda(s)\rightarrow
\lambda_{\pm}$ as $s\rightarrow \pm \infty$, where the limiting
behaviour is asymptotically constant (i.e. $d\Lambda/ds\rightarrow 0$
as $s\rightarrow \pm\infty$).  More precisely, 
\begin{asmp}
We assume
$$
\Lambda(s)\in\cP(\lambda_-,\lambda_+)
$$
where we consider the subset of $C^2(\R,(\lambda_-,\lambda_+))$:
$$
\cP(\lambda_-,\lambda_+)=\left\{\Lambda(s)~:~
  \lambda_-<\Lambda(s)<\lambda_+, \; \lim_{s\rightarrow \pm \infty} \Lambda(s)= \lambda_{\pm} \mbox{ and } \lim_{s\rightarrow \pm \infty} \frac{d\Lambda}{ds} = 0\right\}.
$$
\end{asmp}
Figure~\ref{fig:lambdas} gives two examples of parameter shifts that we consider, though clearly a parameter shift may have an arbitrary number of maxima and minima.

The {\em set of equilibria} for the bifurcation diagram of (\ref{eq:odeaut}) is
$$
\cX=\{(x,\lambda)~:~f(x,\lambda)=0~\mbox{ and }~\lambda_-\leq \lambda\leq \lambda_+\}\subset \R^n\times \R.
$$
Writing $df(x,\lambda)$ to denote the linearization of $f$ with
respect to $x$, we identify the subsets of $\cX$ that are linearly
{\em stable}, {\em unstable} and {\em bifurcation equilibria}:
$$
\begin{array}{l}
\cX_{\stab}=\{(x,\lambda)\in \cX~:~\max\{\re(\sigma(df(x,\lambda)))\}<0\},\\
\cX_{\unst}=\{(x,\lambda)\in \cX~:~\max\{\re(\sigma(df(x,\lambda)))\}>0\},\\
\cX_{\bif}=\{(x,\lambda)\in \cX~:~\max\{\re(\sigma(df(x,\lambda)))\}=0\}
\end{array}
$$
where $\sigma(M)$ denotes the spectrum of the matrix $M$.
Note that, in addition to generic bifurcations, $\cX_{\bif}$ may
include degenerate bifurcation points and even ``bifurcations'' that
do not involve any topological change in the phase portrait, e.g.
$\frac{d}{dt} x=-\lambda^2 x-x^3$ at $\lambda=0$.

We say $\cX$ is {\em regular} if (i) the  set of bifurcation points are
isolated ($\cX_{\bif}$ has no accumulation point) and (ii) the only
asymptotically stable recurrent sets (attractors) are equilibria.
The first property of a regular $\cX$ is implied by the Kupka-Smale
theorem for flows \cite[Thm 7.2.13]{Katok-Hasselblatt} under generic
assumptions, though our assumption is less restrictive and allows the
possibility of high codimension bifurcations in $\cX_{\bif}$. The
second property of a regular $\cX$ holds for one-dimensional systems,
but does not hold in general for higher dimensional systems.

\begin{asmp}
We assume $\cX$ is regular. Furthermore, for ease of exposition, we assume
w.l.o.g. that $\lambda_-<\lambda_+$ throughout and that there are no
bifurcations at $\lambda_{\pm}$.
\end{asmp}

Henceforth, we will denote the closure of $U$ by $\overline{U}$ and define:

\begin{defn}
  We say a $C^2$-smooth curve $X(\lambda)$ in $\cX$ is a {\em branch}
  if it does not intersect with $\cX_{\bif}$ except possibly at end
  points. We say it is a {\em stable branch} if it is contained in
  $\overline{\cX}_{\stab}$.
\end{defn}

\begin{defn}
Given a parameter shift $\Lambda(s)\in \cP(\lambda_-,\lambda_+)$, we
say a continuous curve $(X(s),\Lambda(s))$ that limits to some $(X_{\pm},\lambda_{\pm})$ as $s\rightarrow \pm\infty$ and whose image lies within $\cX$ is a {\em path}. We say it is a {\em stable path} if its image lies within
  $\overline{\cX}_{\stab}$
\end{defn}

Note that branches and paths are defined in terms of solutions to the
autonomous system (\ref{eq:odeaut}) and, to some extent, the shape of
the parameter shift $\Lambda(s)$, but are independent of the rate $r$.
A path can traverse several branches that meet at isolated points in
$\cX_{\bif}$, where the path may be nonsmooth. It may also traverse the same branch several times.  
As illustrated in Figure~\ref{fig:lambdas}, some $\Lambda(s)$ can be monotonic functions of $s$ and
visit each point on a branch at just one value of $s$, others may have internal minima
or maxima and visit the same point on a branch at more than one
value of $s$. We mention a minor result that we use in a later theorem.

\begin{lemma}
If $\cX$ is regular then a stable path traverses a finite number of smooth stable branches.
\end{lemma}

\proof Regularity of $\cX$ in particular means that there are no
accumulations of bifurcation points $\cX_{\bif}$. Suppose that a
stable path $(X(s),\Lambda(s))$ traverses an infinite number of smooth
stable branches and hence passes through an infinite number of bifurcation
points. As the path is continuous and bifurcation points are isolated
on the path, there must be an accumulation as $s\rightarrow \infty$ or
$-\infty$. However, as $X$ has well-defined end-points in these
  limits $X_\pm = \lim_{s\rightarrow\pm\infty}X(s)$, one of
$X_{\pm}$ is an accumulation point, which gives a contradiction.  \qed

\subsection{Local pullback attractors for the non-autonomous system}

We are interested in the influence of parameter shifts
$\Lambda(s)\in\cP(\lambda_-,\lambda_+)$ and rates $r$ on
pullback attractors of the non-autonomous system (\ref{eq:ode}).
Note that the limit $r\rightarrow 0$
is singular in the sense that $r=0$ does not exhibit any change in
parameter values with time. In terms of the cocycle $\Phi_{s,t}$, the following
definition is a local version of the pullback point attractor
discussed in \cite[Defn 3.48(ii)]{KloedenRasmussen} that we have adapted to the context of parameter shifts:

\begin{defn}
We say a solution $x(t)$ of (\ref{eq:ode}) is a (local) {\em pullback
  attractor} if there is a bounded open set $U\subset \R^n$ and a time $T>0$ such that (i) $x(t)\in U$ for all $t<-T$, and (ii) for any $y\in U$ and $t\in\R$ we have
\begin{equation}
|\Phi_{t,t-s}(y)-x(t)|\rightarrow 0 \mbox{ as }s\rightarrow \infty.
\label{eq:pullback}
\end{equation}
\end{defn}

In contrast to \cite[Defn 3.48(ii)]{KloedenRasmussen} we do not require that the pullback attractor attracts a non-autonomous open set $U(t)$.
For the asymptotically autonomous case we consider here, we can relate each attractor for the past limit autonomous system to a pullback attactor for the non-autonomous system via the following theorem, whose proof is given in
Appendix~A.\footnote{We are indebted to an anonymous referee for suggesting this method of proof.}

\begin{theorem}
\label{thm:pullback}
For any $X_-$ such that $(X_{-},\lambda_{-})$ is in $\cX_{\stab}$,
there is a unique trajectory $\tilde{x}_{pb}(t)$ of (\ref{eq:ode}) with $\tilde{x}_{pb}(t)\rightarrow X_{-}$ as $t\rightarrow
-\infty$, and this trajectory is a (local) pullback point attractor
for (\ref{eq:ode}). Moreover there is a neighbourhood $U$ of $X_-$ and $T>0$ such that $\tilde{x}(t)$ is the only trajectory remaining in $U$ for all $t< -T$.
\end{theorem}

Theorem~\ref{thm:pullback} shows that any trajectory that limits to a
linearly stable equilibrium in the past is a pullback
point attractor, and these pullback attractors are in one-to-one
correspondence with choices of $X_-$ such that $(X_-,\lambda_-)\in
\cX_{\stab}$. Note that this result is similar to \cite[Thm
3.54(i)]{KloedenRasmussen}, except that their
assumption that $x=0$ is invariant means that $\tilde{x}_{pb}(t)=0$. 
Any small enough open neighbourhood $U$ of $X_-$ can be
used to define the corresponding pullback attractor of (\ref{eq:ode}),
in particular by choosing small enough $\epsilon>0$, the unique
pullback attractor with past limit $X_-$ can be expressed as
\begin{equation}
\tilde{x}_{pb}^{[\Lambda,r,X_-]}(t):= \bigcap_{s>0}\bigcup_{u>s} \Phi_{t,t-u}(B_{\epsilon}(X_-))
\label{eq:pb}
\end{equation}
for a given parameter shift $\Lambda(s)\in\cP(\lambda_-,\lambda_+)$,
rate $r>0$ and initial state $X_-$, where $B_\epsilon$ denotes a $n-$dimensional ball of radius $\epsilon$.
We do not use (\ref{eq:pb}) in the
following, but note that it may be useful as a basis for numerical
approximation of the pullback attractor that limits to a given $X_-$ for
$t\rightarrow -\infty$.

\subsection{Tracking of stable branches by pullback attractors}

We define two notions of tracking of branches of stable equilibria of
(\ref{eq:odeaut}) by attractors of (\ref{eq:ode}). The strongest notion
of tracking we consider is {\em $\epsilon$-close tracking}:

\begin{defn}
  For a given $\epsilon>0$, we say a (piecewise) solution
  $x(t)$ of (\ref{eq:ode}) for some $r>0$ and $\Lambda(s)\in\cP$ {\em
    $\epsilon$-close tracks} the stable path $(X(s),\Lambda(s))$ if
\begin{equation}
|x(t)-X(r t)|<\epsilon
\label{eq:closetracking}
\end{equation}
for all $t\in\R$.
\end{defn}

The following is another, somewhat weaker notion of tracking, which we call {\em end-point tracking}:

\begin{defn}
Consider a stable path $(X(s),\Lambda(s))$ from $(X_-,\lambda_-)$ to $(X_+,\lambda_+)$. We
say a solution $x(t)$ of (\ref{eq:ode}) {\em end-point tracks} this path if it  satisfies
$$
\lim_{t\rightarrow \pm \infty}x(t)=X_{\pm}.
$$
\end{defn}
Note that end-point tracking implies that for all $\epsilon>0$ there is a $T>0$ such that (\ref{eq:closetracking}) holds for all $|t|>T$. The following result gives a sufficient condition for close- and end-point tracking of a stable branch.

\begin{lemma}
	\label{lem:pbtrack}
	Suppose that $\Lambda(s)\in \cP(\lambda_-,\lambda_+)$ and
        $(X(s),\Lambda(s))$ is a stable path that traverses 
        a stable branch from $(X_-,\lambda_-)$ to
        $(X_+,\lambda_+)$ bounded away from $\cX_{\bif}$. Then
        for any $\epsilon>0$ there is a $r_0>0$ such that for all
        $0<r<r_0$ the pullback attractor $x_{pb}^{[\Lambda,r,X_-]}(t)$
        satisfies $\epsilon$-close tracking of $(X(s),\Lambda(s))$. The pullback attractor 
        also end-point tracks the stable path $(X(s),\Lambda(s))$ for all sufficiently small $0<r$.
\end{lemma}

\proof Let $(X(s),\Lambda(s))$ be the stable path whose
  image is contained within the stable branch $X(\lambda)$ bounded
  away from $\cX_{\bif}$ and we write $X(s)=X(\Lambda(s))$.
% Because the branch is bounded away from $\cX_{\bif}$ there is a
% $\rho>0$ such that
% $\max\{\re(\sigma(df(X(s),\Lambda(s))))\}<-\rho<0$ for all $s\in\R$.
 Now consider the augmented system on $\R^{n+1}$ given by
\begin{equation}
\label{eq:augmented}
\begin{split}
\frac{dx}{dt} & = f(x,\Lambda(s))\\
\frac{ds}{dt} & = r.
\end{split}
\end{equation}
For $r=0$ this has a compact invariant critical manifold
$$
M_{0}:=\{(X(\Lambda(s)),s)~:~s\in\R \}
$$
that is foliated with (neutrally stable) equilibria. Note that the section of $M_0$ with constant $s$ limits to $X_\pm$ for $s\rightarrow\pm\infty$. As the branch $X(\lambda)$ is
bounded away from any bifurcation points, there is a $\rho>0$ such
that the normal Jacobian satisfies
$\max\{\re(\sigma(df(X(s),\Lambda(s))))\}<-\rho<0$ for all $s\in \R$. In consequence, the one-dimensional invariant manifold $M_{0}$ is uniformly normally hyperbolic. One can check that $M_{0}$ is a smooth Riemannian manifold of bounded geometry. Applying\footnote{We use $k=2$ and $\alpha=0$ in the notation of \cite{Eldering} and note that because we restrict to $\cX_{\stab}$ there is an empty unstable bundle.} \cite[Theorem 3.1]{Eldering}, there is an $r_0(\epsilon)>0$ such that for any $r<r_0$ the system (\ref{eq:augmented}) has a unique normally hyperbolic invariant manifold $M_{r}$, where the $C^1$ distance from $M_{0}$ to $M_{r}$ does not exceed  $\epsilon$.

Moreover, this manifold contains a single trajectory parametrized by $s\in\R$.  In summary, for any $\epsilon>0$ and for all $0<r<r_0(\epsilon)$ there is a trajectory $(\tilde{x}(t),rt)\in M_r$ of (\ref{eq:augmented}) with $|\tilde{x}(t)-X(\Lambda(rt))|<\epsilon$ for all $rt\in\R$.

By Theorem~\ref{thm:pullback} there is a neighbourhood $U$ containing $X_-$ such that only one trajectory $\tilde{x}(t)$ satisfies $\tilde{x}(t)\in U$ as $t\rightarrow -\infty$. Fixing $\epsilon>0$ and $T>0$ such that $\tilde{x}(t)\in B_{\epsilon}(X(t))\subset U$ for all $t<-T$ we can see that the trajectory obtained is the pullback attractor
$$
\tilde{x}(t)=x_{pb}^{[\Lambda,r,X_-]}(t).
$$
In particular $\tilde{x}(t)\rightarrow X_-$ as $t\rightarrow -\infty$. 

Now consider the behaviour of $\tilde{x}(t)$ in the limit $t\rightarrow \infty$. If we set $\tilde{x}(t):=X_+ +z(t)$ then
$$
\frac{dz}{dt}=f(X_++z,\lambda_+)+G(z,t)
$$
where $G(z,t)=f(X_++z,\Lambda(rt))-f(X_++z,\lambda_+)$. Noting that 
\begin{equation}
\frac{dz}{dt}=df(X_+,\lambda_+)z+G(z,t)+O(|z|^2)
\label{eq:dzdt}
\end{equation}
where the error term is uniform in $t$. For any $\epsilon>0$ define
$$
g(t):=\sup\{|G(z,t')|~:~|z|<\epsilon ~\mbox{ and }~ t'>t\}
$$
and continuity of $f$ implies that $g(t)\rightarrow 0$ as $t\rightarrow \infty$. For any $\epsilon>0$,
$|\tilde{x}(t)-X(\Lambda(rt))|<\epsilon/2$ for all $r<r_0(\epsilon/2)$
and $t$. The fact that $\lim_{s\rightarrow \infty}X(\Lambda(s))=X^+$ means there is a $T(\epsilon)>0$ such that
$|X_+-X(\Lambda(rt))|<\epsilon/2$ for $t>T(\epsilon)$. Hence
$|z(t)|<\epsilon$ for all $t>T(\epsilon)$. Finally, note
that from the linear stability of $\dot{z}=df(X_+,\lambda_+)z$, there
is an $\epsilon>0$ such that if $|z(T)|<\epsilon$ then
$z(t)\rightarrow 0$ as $t\rightarrow \infty$. It follows that
$\tilde{x}(t)\rightarrow X_+$ as $t\rightarrow \infty$ and so the
solution $\tilde{x}(t)$ end-point tracks the branch.  \qed

\subsection{Tracking of stable paths by pseudo-orbits}

One might wish to characterize a wider class of stable paths that can potentially
be tracked by including stable paths that  pass through bifurcation points where 
stable branches meet. To do this we have to widen our notion of tracking to consider
tracking by pseudo-orbits.

Let us consider the possible stable paths that may be accessible for a
given $\Lambda$.

\begin{defn}
	Given a parameter shift $\Lambda(s)\in \cP(\lambda_-,\lambda_+)$, we
	say points $(X_-,\lambda_-)$ and $(X_+,\lambda_+)$ on $\cX_{\stab}$
	are $\Lambda$-{\em connected} if they lie on the same stable path
	$(X(s),\Lambda(s)).$%\subset
                                %\overline{\cX}_{\stab}$.
\end{defn}

Note that $(X_-,\lambda_-)$ and $(X_+,\lambda_+)$ are
$\Lambda(t)$-connected if and only if they are $\Lambda(r
t)$-connected for any $r>0$, so the definition does not depend on $r$. For given $(X_-,\lambda_-)$ and $(X_+,\lambda_+)$ that are $\Lambda$-connected, it is not
necessarily the case that there will be tracking by orbits of the non-autonomous system for some $r>0$. In particular, there is no guarantee that the pullback attractor $x_{pb}^{[\Lambda,r,X_-]}(t)$ from Theorem~\ref{thm:pullback} remains close to $X(rt)$, that
$$
\lim_{t\rightarrow\infty}x_{pb}^{[\Lambda,r,X_-]}(t)=X_+
$$
or even that $\lim_{t\rightarrow\infty}x_{pb}^{[\Lambda,r,X_-]}(t)$ is an 
attractor for the future limit system. Indeed the example
in Section~\ref{sec:lambda} shows that such a trajectory
may not exist: achieving the appropriate branch switching may
not be possible even on varying $r$ and $\Lambda$.

As discussed later in Section~\ref{sec:lambda} there may be several
stable paths from $X_-$ that end at different $X_+$. We show that
there will be tracking of any stable path by a {\em pseudo-orbit},
namely an orbit with occasional arbitrarily small adjustments. Recall
the following definition of a pseudo-orbit (N.B. there are several
possible definitions of pseudo-orbit for flows; e.g. see
\cite{Palmer}). We say $x(t)$ is an $(\epsilon,T)$-pseudo-orbit
\cite{Shub} of (\ref{eq:ode}) on $[T_a,T_b]$ if there is a finite
series $t_i\in [T_a,T_b]$ for $i=1,\ldots,M$, with $t_1=T_a$ and
$t_M=T_b$, such that for all $i=1,\ldots,M-1$,
\begin{itemize}
\item $t_{i+1}>t_{i}+T$
\item $x(t)$ is a trajectory of (\ref{eq:ode}) on the time-interval $[t_{i},t_{i+1})$
\item $|x(t_{i+1})-\lim_{t\rightarrow t_{i+1} -} x(t)|<\epsilon$.
\end{itemize}
We say $x(t)$ defined for $t$ in some (possibly infinite)
interval is an $(\epsilon,T)$-pseudo-orbit if it is an
$(\epsilon,T)$-pseudo-orbit on all finite sub-intervals. For convenience we refer to an $(\epsilon,1)$-pseudo-orbit as an $\epsilon$-pseudo-orbit.

\begin{theorem}
\label{thm:adiabatic}
Suppose that $(X_-,\lambda_-)$ and $(X_+,\lambda_+)$ are
$\Lambda$-connected by a stable path $(X(s),\Lambda(s))$ for some $\Lambda(s)\in \cP(\lambda_-,\lambda_+)$. Then
for any $\epsilon>0$ there is a $r_0>0$ such that for all
$r<r_0$ there is an $\epsilon$-pseudo-orbit $\tilde{x}(t)$ of (\ref{eq:ode}) with
\begin{equation}
|\tilde{x}(t)-X(r t)|<\epsilon
\label{eq:tracking}
\end{equation}
for all $t\in\R$.
\end{theorem}

This Theorem (proved below) suggests that if we consider the system perturbed by noise of
arbitrarily low amplitude, for any noise amplitude there will be
an initial condition, a rate and a realization of additive noise such that
the perturbed system has a trajectory that end-point tracks the
given stable path. Under additional assumptions, we believe that stable paths will
be $\epsilon$-close tracked with positive probability. Before proving Theorem~\ref{thm:adiabatic},
we give a special case whose proof uses similar argument to that of Lemma~\ref{lem:pbtrack}.

\begin{lemma}
\label{lem:adiabatic}
Suppose that $\Lambda(s)\in \cP(\lambda_-,\lambda_+)$ and $(X(\Lambda(s)),\Lambda(s))$ traverses a
single stable branch that is bounded away from
$\cX_{\bif}$ in the interval $s\in[s_a,s_b]$ (we include the
possibility $s_a=-\infty$ or $s_b=\infty$). Then for any $\epsilon>0$
there is a $r_0>0$ such that for all $0<r<r_0$ there is a trajectory
$\tilde{x}(t)$ that satisfies $\epsilon$-close tracking of this branch in the sense
of (\ref{eq:tracking}) for $t\in[s_a/r,s_b/r]$.
\end{lemma}

\proof This is similar to that of Lemma~\ref{lem:pbtrack} except that instead of $s\in\R$
we consider $s\in[s_a,s_b]$ and $\tilde{x}$ is not necessarily a pullback attractor
for the system, or even unique. Consider the case of $s_a$ and $s_b$ finite and note
that if the closed set $\{(X(s),\Lambda(s)):s\in[s_a,s_b]\}$
does not intersect the closed set $\cX_{bif}$ then they must be isolated by
a neighbourhood. Hence there is
a $\rho>0$ such that $\max\{\re(\sigma(df(X(s),\Lambda(s))))\}<-\rho<0$ for all $s\in S:=[s_a,s_b]$. By considering an
augmented system as in the proof of Lemma~\ref{lem:pbtrack} we can apply Fenichel's theorem \cite{FEN79} 
to conclude that for small enough $r$
there is a trajectory $\tilde{x}(t)$ that $\epsilon$-tracks the branch for $rt=s\in S$.
\qed

~

\proof {\bf (of Theorem~\ref{thm:adiabatic})} 
Let $(X(s),\Lambda(s))$ be the stable path that starts at
$(X_-,\lambda_-)$ and ends at $(X_+,\lambda_+)$. Because $\cX$ is
regular, we can assume that there are at most finitely many
bifurcation points on the path $(X(s),\Lambda(s))$.
However, it is possible that $\Lambda(s)$ may remain at bifurcation points for a set of $s$ that may include intervals: there is a finite set of
points $\{s^{\pm}_1,\cdots,s^{\pm}_k\}$ with
$$
s^-_1\leq s^+_1\cdots<s^-_k\leq s^+_k
$$
such that $\Lambda(s)$ is constant and
$(X(s),\Lambda(s))\in\cX_{\bif}$ when $s\in [s^-_i,s^+_i]$
for some $i$. Note that the generic (transversal) case is $s^-_i=s^+_i$ but we allow a more general case where the parameter can ``linger'' at bifurcation points for a time.

We proceed to show that for any $\epsilon>0$ there is an $r_0(\epsilon)$ such that it is possible to construct an $\epsilon$-pseudo-orbit $\tilde{x}(t)$ that
shadows the stable path $X(r t)$ for all $r<r_0$. All
discontinuities of the constructed pseudo-orbit will be near $t\in 
[s^-_j/r,s^+_j/r]$. We show this in detail for the finite interval
$s\in[s^+_{j-1},s^-_j]$: similar proofs hold for the semi-infinite
intervals $(-\infty,s^-_1]$ and $[s^+_k,\infty)$ by using uniformity of 
linear stability near the endpoints $X_{\pm}$. Let us define
  $t^{\pm}_j:=s^{\pm}_j/r$ and
  $(X_j,\lambda_j):=(X(s^{\pm}_j),\Lambda(s^{\pm}_j))$.

Pick any $\epsilon>0$. Continuity of $X(s)$ means we can find a
$\delta$ such that
\begin{equation}
\label{eq:ineq0}
|X(s)-X_j|<\frac{\epsilon}{4}.
\end{equation}
for all $s\in[s^-_j-\delta,s^-_j]$. 

Pick any $\eta$ with $0<\eta<\delta$; Lemma~\ref{lem:adiabatic}
implies that there is a rate $r_0>0$ (depending on $\epsilon$ and
$\eta$) and a trajectory $x(t)$ such that for all $r<r_0$ there is a
trajectory $x_j(t)$ of (\ref{eq:ode}) satisfying
\begin{equation}
\label{eq:ineq1}
|x(t)-X(r t)|<\frac{\epsilon}{4}
\end{equation}
on the interval $r t\in[s^+_{j-1}+\eta,s^-_{j}-\eta]$.

The fact that $(X_{j},\lambda_{j})$ is a bifurcation point
($f(X_j,\lambda_j)=df(X_j,\lambda_j)=0$) and $f\in C^2$ means that, by Taylor expansion, for any
$\xi,\zeta>0$ there are $M,N>0$ such that
$$
|f(x,\lambda)|<M|x-X_{j}|^2+N|\lambda-\lambda_{j}|
$$
whenever $|x-X_{j}|<\xi$ and
$|\lambda-\lambda_{j}|<\zeta$. Similarly, $\Lambda(s)\in C^1$ means
there is a $P>0$ such that
$$
|\Lambda(s)-\lambda_j|<P d(s,[s^-_j,s^+_j])
$$
for all $s\in [s^-_j-\eta,s^+_j+\eta]$, where
$d(s,A)=\inf\{|s-t|~:~t\in A\}$. Choosing $\zeta=P\eta$ we have
\begin{equation}
|f(x,\Lambda(r t))|< M\xi^2+NP\eta
\end{equation}
for all points in
\begin{equation}
\label{eq:nbhd}
I_{j}(\xi,\eta):=\{(x,s)~:~|x-X_j|<\xi\mbox{ and }s\in[s^-_j-\eta,s^+_j+\eta]\}.
\end{equation}

Considering a forward trajectory $x(t)$ of (\ref{eq:ode}) we have 
\begin{equation}
\label{eq:ineqa}
|x(t)-x(t')|< (M\xi^2+NP\eta) |t-t'|
\end{equation}
for $t>t'$, if $(x(t),rt)\in I_{j}(\xi,\eta)$ between $t$ and $t'$. Hence by choosing $\xi$ such that
$$
\epsilon<\xi <\sqrt{\frac{\epsilon}{16M}} \mbox{ and }\eta<\frac{\epsilon}{16NP} 
$$
(which is possible as long as $\epsilon<\frac{1}{16M}$) we can ensure that:
\begin{itemize}
\item There is containment $I_{j}(\xi,\epsilon)\subset I_{j}(\xi,\eta)$.
\item From (\ref{eq:ineqa}), the inequality
\begin{equation}
\label{eq:ineq2}
|x(t)-x(u)|<\frac{\epsilon}{4} \mbox{ for all } t\in[u,u+2]
\end{equation}
holds for any trajectory $x(t)$ such that $|x(t)-X_j|<\epsilon$ for all $t\in[u,u+2]$ where $(x(t),rt)\in I_j(\xi,\eta)$. 
\end{itemize}
Because of (\ref{eq:ineq2}), when $rt\in[s^-_j-\eta,s^+_j+\eta]$ we can find a finite number of trajectory segments, each of time
length $T\geq 1$, starting at $x(u)=X(ru)$, such that $|x(t)-x(u)|<\epsilon/4$ for all $t\in[u,u+2]$.

\begin{figure}%
\centerline{\includegraphics[width=8cm]{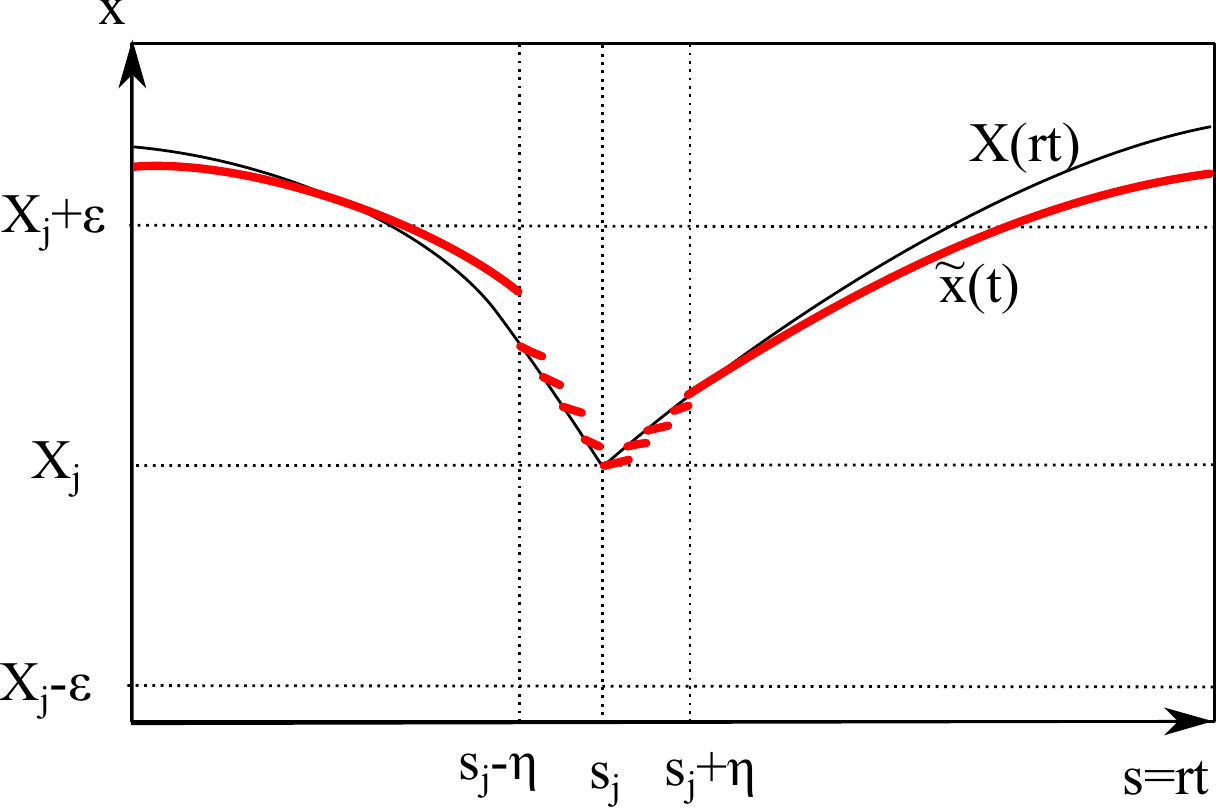}}
\caption{ Construction of a pseudo-orbit $\tilde{x}(t)$ (red) that $\epsilon$-close tracks a
  stable path $X(r t)$ (black) through a bifurcation point. We
  sketch the case where $s_j=s^-_j=s_j^+$ (i.e. the path does not
  linger at the bifurcation point). The short
  sections near the bifurcation point have duration between $1$ and
  $2$ time-units so that the width of the smaller intervals of continuity in $s$ lie
  between $1/r$ and $2/r$. If $s_j^-<s_j^+$ then we choose a pseudo-orbit that remains at $X_j$ over this interval.
}%
\label{fig:pseudo}
\end{figure}

Hence there is a choice of $\eta>0$ and $r_j>0$ such that, if $r<r_j$
then there is a pseudo-orbit $\tilde{x}(t)$ with the desired property on
the interval $rt\in [s^-_{j}-\eta,s^-_{j+1}-\eta]$. More precisely, for any $rt\in S:=[s^-_{j}-\eta,s^+_{j}+\eta]$ we can find an increasing sequence of $s_{j,k}$ and a piecewise solution $\tilde{x}(t)$ with $s_{j,1}=s^-_{j}-\eta$, $s_{j,M}=s^+_{j}+\eta$ such that $s_{j,k+1}>s_{j,k}+1$. Combining the estimates (\ref{eq:ineq0},\ref{eq:ineq1}) and (\ref{eq:ineq2}) we have for this solution
\begin{eqnarray}
|\tilde{x}(t)-X(r t)| & <& |\tilde{x}(t)-\tilde{x}(u)|+|\tilde{x}(u)-X(r u)|+|X(r u)-X_j|+|X_j-X(r t)|\nonumber\\
&<& \frac{\epsilon}{4}+\frac{\epsilon}{4}+\frac{\epsilon}{4}+\frac{\epsilon}{4}  \nonumber\\
& = & \epsilon.\label{eq:xestimate}
\end{eqnarray}
for all $rt\in [s^-_{j}-\eta,s^-_{j+1}-\eta]$ except where $rt=s_{j,k}$ for some $k$.

On choosing $r<\min\{r_j~:~j=1,\ldots,k\}$ we have a piecewise defined path $\tilde{x}(t)$ which either satisfies (\ref{eq:ode}) at $t$ and (\ref{eq:xestimate}) holds, or
\begin{eqnarray*}
|\lim_{u\rightarrow t-}\tilde{x}(u) -\lim_{u\rightarrow t+}\tilde{x}(u)| & < & |\lim_{u\rightarrow t-}\tilde{x}(u)-X(r t)| < \frac{\epsilon}{4}<\epsilon.
\end{eqnarray*}
By taking the upper limit at any discontinuity, we note that $\tilde{x}(t)$ is an $\epsilon$-pseudo-orbit as required.
\qed

\subsection{An example: dependence on $\Lambda$-connectedness}
\label{sec:lambda}

\begin{figure}%
\centerline{\includegraphics[width=6cm]{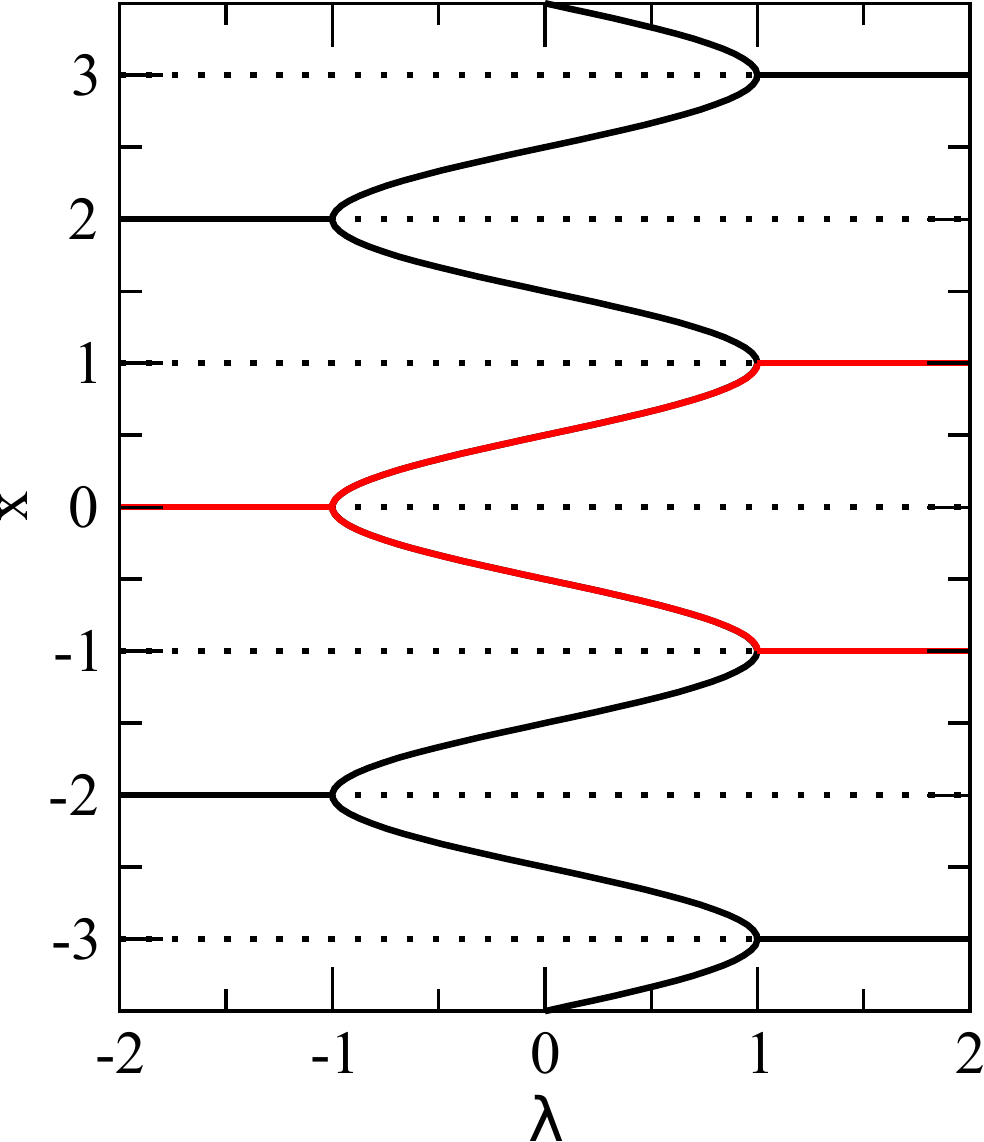}}
\caption{ Bifurcation diagram for (\ref{eq:changeover}); the solid
  lines show the stable branches while the dashed lines are
  unstable. Observe that the only
  possible tracking by a pseudo-orbit that starts at $(X_-,\lambda_-)=(0,-2)$ for a monotonic
  $\Lambda(s)$ as in Figure~\ref{fig:lambdas}(a) will be the red lines shown that finish at
  $(X_+,\lambda_+)=(\pm 1,2)$. However, for $\Lambda(s)$ as in Figure~\ref{fig:lambdas}(b) one can also track
  paths by pseudo-orbits that finish at $(X_+,\lambda_+)=(\pm 3,2)$.
}%
\label{fig:changeover}
\end{figure}

Let us consider the system
\begin{equation}
\frac{dx}{dt}=\sin(x\pi)(\Lambda(r t)+\cos(x\pi))
\label{eq:changeover}
\end{equation}
for $x$, $r>0$ and $\Lambda(s)\in\cP(-2,2)$ this is a simple example
where we can explore which stable paths can be adiabatically tracked
for a given $\Lambda(s)$. If we note that $\lambda_-=-2$ and
$\lambda_+=+2$ then the bifurcation diagram is as in
Figure~\ref{fig:changeover}.  Although there are branches that
continue through for $x=k\in\Z$, there is an exchange of stability by
a combination of pitchfork bifurcations. For this system we see that
the set of stable branches is connected into one component and so for
any $k,l\in\Z$ there is a stable path from $(2k,-2)$ to
$(2l+1,2)$. However, if $|k-l|$ is large then the path 
  will require $\Lambda(s)$ with a large number of internal maxima
and minima.

Now consider any trajectory $x(t)$ of (\ref{eq:changeover}). As the
lines of equilibria at integer values of $x$ remain flow-invariant
subspaces for any choice of $\Lambda$ and $r$, it follows that
$$
\lfloor x(0)\rfloor \leq x(t)< \lfloor x(0)\rfloor+1~~\mbox{ for all }t,
$$
where $\lfloor x\rfloor$ is the smallest integer $k\leq x$.  This
means that the pullback attractor from $X_-=0$ cannot track a stable
path from $(X_-,\lambda_-)=(0,-2)$ to $(X_+,\lambda_+)$ 
where $X_+=\pm1,\pm3,\ldots$, and $\lambda_+=2$ in
Figure~\ref{fig:changeover}. Indeed, it is easy to see that the
pullback attractor in this case is simply the constant
$\tilde{x}_{pb}^{[\Lambda,r,X_-]}=X_-$ and in all cases this is in
fact a repellor for the future limit system.

However, Theorem~\ref{thm:adiabatic} implies that a {\em pseudo-orbit}
can track such a path - it allows the possibility of arbitrarily small
adjustments along the way. Starting at $(X_-,\lambda_-)=(0,-2)$ and
allowing arbitrarily small adjustments, there are two possible
finishing points for a $\Lambda(s)$ such as
Figure~\ref{fig:lambdas}(a) and four possible finishing points for a
$\Lambda(s)$ such as Figure~\ref{fig:lambdas}(b).

%%%%%%%%%%%%%%%%%%%%%%%%%%%%%%%%%%%%%%%%%%%%%%%%%%%%%%%%%%%%%%%%%%%%
\subsection{Bifurcation-equivalent systems}
\label{sec:equivalent}

Since the discussion so far only depends on branches of equilibria,
their stabilities and bifurcations, many of the properties carry
across to systems that have the same bifurcation diagram. Suppose
that
\begin{equation}
\label{eq:odeequiv}
\frac{dy}{d\tau}=\rho(y) f(y,\Lambda(r \tau))
\end{equation}
where $y\in\R^n$ and $\rho\in C^2(\R^n,\R)$ is a strictly positive
scalar function that is bounded below. 
The associated autonomous system
\begin{equation}
\label{eq:autodeequiv}
\frac{dy}{d\tau}=\rho(y) f(y,\lambda)
\end{equation}
will have the same trajectories in the phase portrait, different
magnitude but the same sign eigenvalues, and hence the same
bifurcation diagram as (\ref{eq:odeaut}). For the same reason systems
(\ref{eq:ode}) and (\ref{eq:odeequiv}) will have the same
$\Lambda$-connected points $X_-$ and $X_+$.  They also have orbits
that are equivalent in the following (rather weak) sense:

\begin{theorem}
\label{thm:bifequiv}
Consider non-autonomous systems (\ref{eq:ode})  and (\ref{eq:odeequiv}).
\begin{itemize}
\item[(a)]
Suppose that for some $\Lambda(rt)\in\cP(\lambda_-,\lambda_+)$ and
$r>0$ system (\ref{eq:ode}) has a solution $x(t)$. Then there exists a
(state-dependent) rescaling of time $T(\tau)$ with $T(0)=0$, and a
$\tilde{y}(\tau) = x(T(\tau))$ that solves the equivalent system
(\ref{eq:odeequiv}) with $\tilde{\Lambda}(r\tau)=\Lambda(rT(\tau))$
replacing $\Lambda(r\tau)$.
\item[(b)] Conversely, suppose that for some
  $\Lambda(r\tau)\in\cP(\lambda_-,\lambda_+)$ and $r>0$ system
  (\ref{eq:odeequiv}) has a solution $y(\tau)$. Then, this $y(\tau)$
  defines a (state-dependent) time rescaling $T(\tau)$ with $T(0)=0$
  such that the equivalent system (\ref{eq:ode}) has a solution
  $\tilde{x}(t)=y(T^{-1}(t))$ with $\tilde{\Lambda}(rt)=\Lambda(
  rT^{-1}(t))$ replacing $\Lambda(rt)$.
\end{itemize}
\end{theorem}

\proof
To show $(a)$, define a (state-dependent) time rescalling $T(\tau)$
through a coupled system
\begin{align*}
\frac{dT}{d\tau} &=\rho(y),\\
\frac{dy}{d\tau} &=\rho(y)\,f\left(y,\Lambda(rT)\right),
\end{align*}
with $T(0)=0$ and $\Lambda(rT(\tau))=\tilde{\Lambda}(r\tau)$. Using
the assumption on $x(t)$ one can
verify that $\tilde{y}(\tau)=x(T(\tau))$ together with a suitable
monotonically increasing function $T(\tau)$ solve the coupled system
above:
\begin{align*}
\frac{d\tilde{y}}{dt} & = 
\frac{dx}{dT}\frac{dT}{d\tau} = \rho(\tilde{y}(\tau))\,f\left(x(T(\tau)),\Lambda(rT(\tau))\right) \\
&=
\rho(\tilde{y})\,f\left(\tilde{y},\Lambda(rT)\right) =
\rho(\tilde{y})\,f(\tilde{y},\tilde{\Lambda}(r\tau)).
\end{align*}
Hence, $\tilde{y}(\tau)$  also solves (\ref{eq:odeequiv})
with $\tilde{\Lambda}(r\tau)=\Lambda(rT(\tau))$ replacing
$\Lambda(r\tau)$.

To show $(b)$, for a given $y(\tau)$ define a state-dependent time
rescaling through the integral
$$
T(\tau)=\int_{u=0}^\tau g(y(u))du,
$$
that satisfies $T(0)=0$ and is a monotonically increasing (hence
invertible) function of
$\tau$ whose first derivative $dT/d\tau = \rho(y(\tau))$ is bounded
from below. We can now verify that $\tilde{x}(t=T(\tau))=y(\tau)$
solves (\ref{eq:ode}) with
$\tilde{\Lambda}(rt)=\Lambda\left(rT^{-1}(t)\right)$ replacing
$\Lambda(rt)$.
On the one hand
$$
\frac{dy}{d\tau} = 
\frac{d\tilde{x}}{d\tau} = 
\frac{d\tilde{x}}{dT}\frac{dT}{d\tau} = 
\rho(y(\tau))\frac{d\tilde{x}}{dT}=
\rho(y(\tau))\frac{d\tilde{x}}{dt}.
$$
On the other hand, the assumption on $y(\tau)$ gives
\begin{align*}
\frac{dy}{d\tau} &= 
\rho(y(\tau))\,f(y(\tau),\Lambda(r\tau))=
\rho(y(\tau))\,f(\tilde{x}(t=T(\tau)),\Lambda(r\tau))\\
&=
\rho(y(\tau))\,f\left(\tilde{x}(t),\Lambda(rT^{-1}(t))\right)\\
&=
\rho(y(\tau))\,f\left(\tilde{x}(t),\tilde{\Lambda}(rt)\right).
\end{align*}
Since $\rho(y)>0$, we have that 
$$
\frac{d\tilde{x}}{dt} = f\left(\tilde{x},\tilde{\Lambda}(rt)\right).
$$
\qed

In particular, this means that a tracking trajectory in an equivalent
system implies that there is a tracking trajectory in
the first system for a different $\Lambda$. Note
  that as the transformation from $\Lambda$ to $\tilde{\Lambda}$ does
  depend on the trajectory chosen,
%and is not always possible to obtain [case (a)],  
  one needs to be cautious about applying Theorem~\ref{thm:bifequiv}.
 We also note that varying $r$ can lead to different parameter shifts
$\tilde{\Lambda}$.  Nonetheless, it does allow extension of some
results to more general bifurcation equivalent systems.

%%%%%%%%%%%%%%%%%%%%%%%%%%%%%%%%%%%%%%%%%%%%%%%%%%%%%%%%%%%%%%%%%%%%
\section{Bifurcation- and rate-induced tipping}
\label{sec:BRtipping}

We now propose some definitions and results on bifurcation (B-) and
rate (R-) induced tipping in the sense of
\cite{AshWieVitCox2012}.  Recall from \cite{AshWieVitCox2012} that B-tipping is associated with a
bifurcation point of the corresponding autonomous system, while
R-tipping appears only for fast enough variation of parameters and is
associated with loss of tracking of a stable branch.

\begin{defn}
\label{def:btipping}
Suppose that $(X_{\pm},\lambda_{\pm})\in\cX_{\stab}$  and fix $\Lambda(s)\in
\cP(\lambda_-,\lambda_+)$. Suppose for all small enough $r>0$  we have 
$$
\lim_{t\rightarrow \infty} \tilde{x}_{pb}^{[\Lambda,r,X_-]}(t) = X_+.
$$
If $X_+$ and $X_-$ are not $\Lambda$-connected then we say there is
  {\em B-tipping} from $X_-$ for this $\Lambda$.
\end{defn}

Definition~\ref{def:btipping} requires that the system
reaches a point on a stable branch that is not accessible by
following a stable path permitted by the chosen $\Lambda$. 
The following result shows that a bifurcation must indeed be responsible for B-tipping. 
%\red{Note: the sentence here seemed out of place and I moved it to the last paragraph}.

\begin{corol}
\label{cor:suffBtip}
Suppose there is B-tipping for some $\Lambda(s)$ starting at $X_-$. Then there is at least one bifurcation point $\lambda_b$ with
$\lambda_-<\lambda_b<\lambda_+$ on the stable branch starting at $X_-$.
\end{corol}

\proof If there is no such bifurcation point, then there must be a
stable branch spanning from $(X_-,\lambda_-)$ to $(X_+,\lambda_+)$.
For any parameter shift $\Lambda(s)\in\cP(\lambda_-,\lambda_+)$,
Lemma~\ref{lem:pbtrack} implies that we track this branch for all
small enough $r>0$ - a contradiction to the assumption that there is
B-tipping.  \qed

We define irreversible R-tipping as follows:

\begin{defn}
Suppose that $(X_-,\lambda_-)\in \cX_{\stab}$ and fix $\Lambda(s)\in
\cP(\lambda_-,\lambda_+)$. We say there is {\em irreversible
  R-tipping} from $X_-$ on $\Lambda$ if there is a $r_0>0$ and an
$X_+$ that is $\Lambda$-connected to $X_-$ such that the system
end-point tracks a stable path  from $X_-$ to $X_+$ for
  $0<r<r_0$ but not for $r=r_0$, i.e.
$$
\lim_{t\rightarrow \infty} \tilde{x}_{pb}^{[\Lambda,r,X_-]}(t) = 
\left\{\begin{array}{cl}
X_+ & ~\mbox{ for }0<r<r_0\\
Y_+ \neq X_+ & ~\mbox{ for }r=r_0
\end{array}\right.
$$
\end{defn}

Although Theorem~\ref{thm:adiabatic} guarantees $\epsilon$-close tracking by pseudo-orbits for small
enough $r$, an increase in $r$ may lead to irreversible R-tipping as above, but then
further increase of $r$ may lead back to tracking again! Indeed the values of
$r$ that give tracking may be interspersed by several windows of rates
that give tipping; see the example in Section~\ref{subsec:rtipexample} and Figure~\ref{fig:bumpratebifs}.

In addition to {\em irreversible R-tipping} defined above, 
there can be {\em transient R-tipping}.  Recall that our definition of
end-point tracking makes no assumption that $|\tilde{x}(t)-X(t)|$ is
small, except in the limits $t\rightarrow \pm\infty$. This means that
the system may end-point track the stable path for all $r>0$, but
depart from the path temporarily when $r>r_0$ to visit a different
state/path at intermediate times~\cite{Wieczorek_etal_2010}. Such
transient R-tipping is not discussed here.

\subsection{An example with rate-induced tipping for intermediate rates}
\label{subsec:rtipexample}
 
\begin{figure}%
\centerline{\includegraphics[width=15cm]{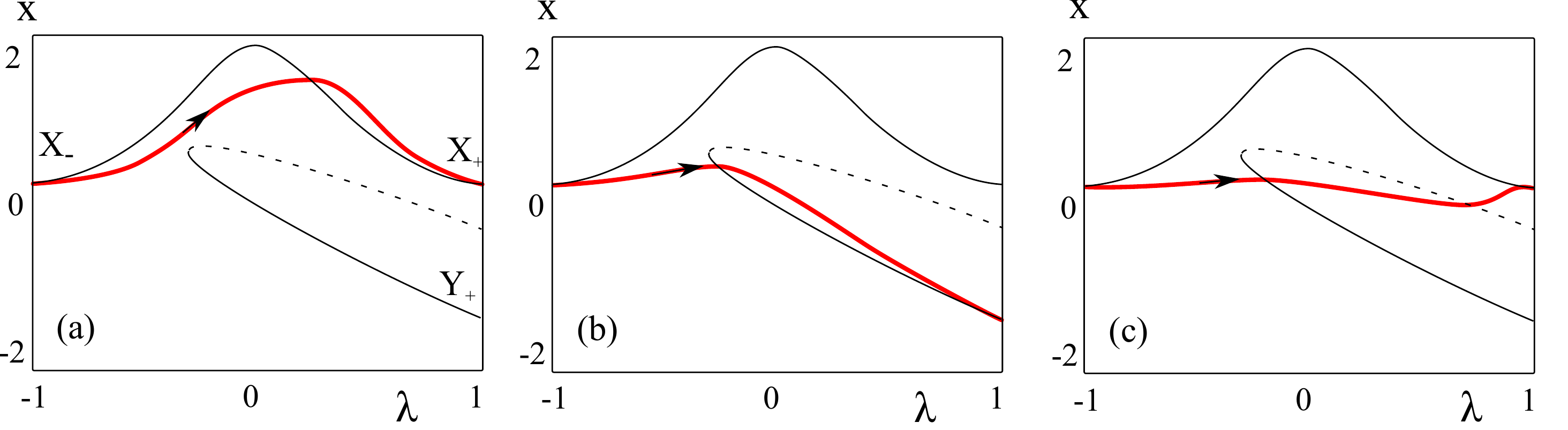}}%
\caption{Bifurcation diagrams showing stable (black solid lines) and unstable (black dashed lines) equilibria for the system (\ref{eq:bump},\ref{eq:bumpfn}) with $\lambda$ between $\lambda_-=-1$ and $\lambda_+=+1$. The red lines show a numerical approximation of the pullback attractor starting at $X_-$ for  parameter shifts (\ref{eq:ppath}) on varying the rate $r$. (a) For $r<r_1$ the pullback attractor end-point tracks the branch through to $X_+$. (b) There is a range of rates $r_1<r<r_2$ where the pullback attractor R-tips to $Y_+$. (c) For $r>r_2$ the pullback attractor again end-point tracks the branch through to $X_+$.}%
\label{fig:bumpratebifs}%
\end{figure}

Consider $x\in \R$ governed by
\begin{equation}
\frac{dx}{dt} = g(x,\Lambda(r t))
\label{eq:bump}
\end{equation}
where $\Lambda(s)\in\cP(\lambda_-,\lambda_+)$ and
\begin{equation}
g(x,\lambda)= -((x+a+b\lambda)^2+c\tanh(\lambda-d))(x-K/\cosh(e\lambda))
\label{eq:bumpfn}
\end{equation}
for $a,b,c,d,e$ and $K$ constants. This example has a particular bifurcation structure as illustrated in Figure~\ref{fig:bumpratebifs} if we choose
$$
a=-0.25,~b=1.2,~c=-0.4,~d=-0.3,~e=3,~K=2.
$$ 
Observe that this is a structurally stable bifurcation diagram and so qualitatively robust to small changes in any parameter, and indeed to any small changes in $g(x,\lambda)$. For $\Lambda(s)\in\cP(\lambda_-,\lambda_+)$ and $r>0$ we take solutions of
\begin{equation}
\frac{d\Lambda}{ds}=-(\Lambda-\lambda_-)(\Lambda-\lambda_+).
\label{eq:ppath}
\end{equation}
Figure~\ref{fig:bumpratebifs} illustrates the appearance of critical rates $0<r_1<r_2$ such that R-tipping only occurs for a finite range $r\in (r_1,r_2)$.

\subsection{Identifying regions of R-tipping in one-dimensional systems}

In this section we assume one state space dimension ($n=1$) throughout and show that geometric properties of the
bifurcation diagram can be used to infer where R-tipping is or is not possible.

Consider a stable path  $(X(s),\Lambda(s))$ for some
$\Lambda(s)\in\cP(\lambda_-,\lambda_+)$. For a fixed
  $v\in\R$, let  $\B(X(v),\Lambda(v))$ be the
basin of attraction of  the stable equilibrium
  $X=X(v)$ for the autonomous
system (\ref{eq:odeaut}) with parameter
  $\lambda=\Lambda(v)$, i.e.
$$
\B(X(v),\Lambda(v)):=\{x_0~:~d(\phi_t(x_0,\lambda),X)\rightarrow 0 \mbox{ as
}t\rightarrow \infty\},
%\B(X,\lambda):=\{x_0~:~d(\phi_t(x_0,\lambda),X)\rightarrow 0 \mbox{ as
%}t\rightarrow \infty\}
$$
where $\phi_t(x_0,\lambda)$ is the trajectory of the autonomous
system (\ref{eq:odeaut}) with  initial condition $x_0$.
We define a notion of {\em forward basin stability} in terms of
$\lambda$-dependent basin of attraction for the attractors of the associated autonomous system on a given stable path.

\begin{defn}
  We say a stable path $(X(s),\Lambda(s))$ for some
  $\Lambda(s)\in\cP(\lambda_-,\lambda_+)$ is {\em forward
    basin-stable} if for every $v\in\R$ the basin of attraction of
  the stable equilibrium $X(v)$ on the path
 contains the closure of the set of all earlier positions of
  this equilibrium along the path:
$$
\overline{\{X(u)~:~u<v\}}\subset \B(X(v),\Lambda(v))\quad\mbox{for
  all}\quad v\in\R.
$$
\end{defn}

Note that forward basin stability: (i) is defined in terms of solutions to the
autonomous system (\ref{eq:odeaut}) and the shape of the parameter shift $\Lambda(s)$, but is independent of the rate $r$, and  (ii) is a property of a stable path rather than a stable branch traversed by the path. We remark that this is a somewhat different notion to basin stability as defined in \cite{Menck2013}: this is a quantity that characterises the relative volume of a basin of attraction and is not a property of a branch of attractors.

In the following, we use the notion of forward basin stability
to give sufficient conditions for irreversible R-tipping in one
dimension to be absent (case 1) or present (cases 2, 3). 

\begin{theorem}
\label{thm:rtipping}
Suppose that $(X_-,\lambda_-)$ and $(X_+,\lambda_+)$ in $\cX_{\stab}$
are $\Lambda$-connected by the stable path $(X(s),\Lambda(s))$ for
some $\Lambda\in \cP(\lambda_-,\lambda_+)$:
\begin{enumerate}
\item If $(X(s),\Lambda(s))$ is forward basin-stable then
$
\tilde{x}_{pb}^{[\Lambda,r,X_-]}(t)
$
end-point tracks $(X(s),\Lambda(s))$ for all $r>0$. Hence in such a case there
can be no irreversible R-tipping from $X_-$ on $\Lambda$.
\item If the stable path $(X(s),\Lambda(s))$ traverses a single stable
  branch $X(\lambda)$, there
is another stable branch\footnote{NB this may be defined only on $[\lambda_0,\lambda_+]$ where $\lambda_-\leq \lambda_0<\lambda_+$.} $Y(\lambda)$ with $Y_+:=Y(\lambda_+) \neq X_+$ and there are $u<v$ such that
$$
X(u) \in \B(Y(\Lambda(v)),\Lambda(v)),
$$
then $(X(s),\Lambda(s))$ is not forward basin-stable and there is a parameter shift
$\tilde{\Lambda}(s)\in\cP(\lambda_-,\lambda_+)$ such that there is
irreversible R-tipping from $X_-$ on this $\tilde{\Lambda}$.
\item 
If the stable path $(X(s),\Lambda(s))$ traverses a single stable branch, and there is a $(Y_+,\lambda_+)\in\cX_{\stab}$ with $X_-$ and $Y_+$
not $\Lambda$-connected such that
$$
X_-\in \B(Y_+,\lambda_+),
$$
then $(X(s),\Lambda(s))$ is not forward basin-stable and there is irreversible R-tipping from $X_-$ for this $\Lambda$.
\end{enumerate}
\end{theorem}

\proof
We define
\begin{align*}
X_{\fut}(u) &:= \overline{\bigcup_{v>u}X(v)}= \left[\inf_{v>u}X(v),\sup_{v>u}X(v)\right],\\
\B_{\fut}(u)& :=\bigcap_{v>u} \B(X(v),\Lambda(v)).
\end{align*}
and note that from these definitions, $X_{\fut}(u) \supset X_{\fut}(v)$ and  $\B_{\fut}(u) \subset \B_{\fut}(v)$ for all $u<v$.

For case 1, let us fix any $r>0$ and suppose that $(X(s),\Lambda(s))$ is forward basin-stable. This implies that $X(u)\in X_{\fut}(u) \cap \B_{\fut}(u)$ and also that $X_-\in \B_{\fut}(u)$ for all $u\in\R$. For any $t_0$ and $x_0\in \B_{\fut}(rt_0)$ we define $x(t):=\Phi_{t,t_0}x_0$ and claim that
$$
\Delta(t):=d(x(t),X_{\fut}(rt_0))
$$
satisfies $\Delta(t)\rightarrow 0$ as $t\rightarrow \infty$. To see this, note that if $\Delta(t)>0$ for any $t>t_0$ then either 
$$
x(t)<\inf_{v>rt_0}X(v) \mbox{  or }x(t)> \sup_{v>rt_0}X(v).
$$
In either case $\frac{d}{dt}{\Delta}(t)<0$ and so $\Delta(t)$ is monotonic decreasing with $t$. If it converges to some $\Delta_0>0$ then we can infer there is another equilibrium (not equal to $X_+$) within $\B(X_+,\lambda_+)$ for the future limit system: this is a contradiction. Hence for any trajectory $x(t)$ such that $x(t)\rightarrow X_-$ as $t\rightarrow \infty$ we have 
$$
d(\Phi_{t,t_0}x,X_{\fut}(rt_0))\rightarrow 0.
$$
Since this holds for all $t_0$ and $\bigcap_{u>0} X_{\fut}(u)=X_+$ we have $x(t)\rightarrow X_+$ as $t\rightarrow \infty$. By Theorem~\ref{thm:pullback} this trajectory is the pullback attractor and hence the result.

For case 2 the existence of such $u,v$ implies that $X(u)\not\in
\B(X(v),\Lambda(v))$ and hence $(X(s),\Lambda(s))$ is not forward basin
stable. We construct a reparametrization
$$
\tilde{\Lambda}(s):= \Lambda(\sigma(s))
$$
using a smooth monotonic increasing $\sigma\in C^2(\R,\R)$ that
increases rapidly from $\sigma(s)=u$ to $\sigma(s)=v$ but increases slowly
otherwise.
More precisely, for any $M>0$ and $\eta>0$ we choose a smooth monotonic function $\sigma(s)$ such that 
$$
\begin{array}{cl}
\sigma(s)=s & \mbox{ for }~~ s<u,\\
1\leq \frac{d}{ds}\sigma(s)\leq M & \mbox{ for }~~ u\leq \sigma(s)\leq v\\ \frac{d}{ds}\sigma(s)=M & \mbox{ for }~~ u+\eta<\sigma(s)<v-\eta,~~\mbox{ and}\\ \frac{d}{ds}\sigma(s)=1 & \mbox{ for }~~ \sigma(s)>v
\end{array}
$$

Using the assumption that the stable path traverses a
single stable branch, we apply Lemma~\ref{lem:pbtrack} to show that for all small enough rates $r>0$ we have
$$
\tilde{x}_{pb}^{[\Lambda,r,X_-]}(u/r)\in \B(Y(v),\Lambda(v)).
$$
Note that by construction,
$$
\tilde{x}_{pb}^{[\tilde{\Lambda},r,X_-]}(t)=\tilde{x}_{pb}^{[\Lambda,r,X_-]}(t)
$$
for any $t\leq u/r$. Indeed, by picking $M$ large and $\eta$ small, we can ensure that 
$$
|\tilde{x}_{pb}^{[\tilde{\Lambda},r,X_-]}(v/r)-Y(\tilde{\Lambda}(v))|
$$
is as small as desired. In particular this means that $\tilde{x}_{pb}^{[\tilde{\Lambda},r,X_-]}(v/r)$ is in the basin of $Y(\tilde{\Lambda}(v))$. Hence $\tilde{x}_{pb}^{[\tilde{\Lambda},r,X_-]}(t)$ will track $Y(\tilde{\Lambda}(t))$ for $t>v/r$ if $r$ is small enough (recall that $s=rt$). 
Applying Lemma~\ref{lem:adiabatic} again for the interval $s\in[v,\infty)$ means that reducing $r$, if necessary, will give a pullback attractor with $\tilde{x}_{pb}^{[\tilde{\Lambda},r,X_-]}(t)\rightarrow Y_+$ as $t\rightarrow\infty$.

For case 3 a similar argument to case 2 implies that the path is not forward basin stable. Given any $\epsilon>0$ there will be a $K>0$ such that if
$$
I:=\overline{\{X(u')~:~u'<-K\}}
$$
then
$$
I\subset [X_--\epsilon,X_-+\epsilon]\subset [X_--2\epsilon,X_-+2\epsilon] \subset
 \bigcap_{v'>K} \B(Y(v'),\Lambda(v')).
$$
By continuity, there is an $M$ such that $|f(x,\lambda)|<M$ for all $x\in I$ and $\lambda\in[\lambda_-,\lambda_+]$. Now pick any 
$$
r> 2\frac{MK}{\epsilon}
$$
so that for the pullback attractor $x(t):=\tilde{x}_{pb}^{[\Lambda,r,X_-]}(t)$ we have $x(-K/r)\in I$ and 
$$
|x(K/r)-x(-K/r)|< 2\frac{MK}{r} <\epsilon.
$$
Therefore 
$$
x(K/r)\in \bigcap_{v'>K} \B(Y(v'),\Lambda(v'))
$$
and so $x(t)\rightarrow Y_+$ as $t\rightarrow \infty$: there is R-tipping on this path for some $0<r_0<r$.
\qed

\subsection{Isolation of B-tipping and R-tipping in one-dimensional systems}
\label{sec:restrictions}

An interesting consequence of Theorem~\ref{thm:rtipping} is that
rate-induced tipping will appear neither for too small a range of
parameters, nor on a segment of a branch that lies too
  close to a saddle-node bifurcation.

\begin{lemma}
\label{lem:rminimal}
If $(X_-,\lambda_-)\in\cX_{\stab}$ then there is a $\nu>0$ such that for $\lambda_+=\lambda_-+\nu$ there is a stable branch $X(\lambda)$ on $(\lambda_-,\lambda_+)$, and no irreversible R-tipping is possible from $X_-$ for any $\Lambda\in\cP(\lambda_-,\lambda_+)$.
\end{lemma}

\proof As $X_-$ is a linearly stable equilibrium, for all small enough $\delta>0$ there will be a $\nu>0$ such that if $\lambda_+=\lambda_-+\nu$ then
$$
B_{\delta}(X_-)\subset \bigcup_{\lambda_-<\lambda<\lambda_+} \B(X(\lambda),\lambda).
$$
Fixing such a $\delta$, continuity of $X(\lambda)$ means that choosing a (possibly smaller) $\nu>0$ we can ensure
$$
\overline{\{X(\lambda)~:~\lambda_-<\lambda<\lambda_+\}} \subset B_{\delta}(X_-)\subset \bigcup_{\lambda_-<\lambda<\lambda_+} \B(X(\lambda),\lambda).
$$
Hence any path on this branch from $\lambda_-$ to $\lambda_+$ will be forward basin stable: by Theorem~\ref{thm:rtipping} case 1, R-tipping is not possible on such a path.
\qed

In fact, a stronger statement may be made near a saddle-node bifurcation.
Suppose there is a saddle-node bifurcation at $(\lambda_0,x_0)$ where the branches exist for $\lambda>\lambda_0$. The unstable and stable branches from the
saddle-node bifurcation will be quadratically tangent to the line
$\lambda=\lambda_0$ in the bifurcation diagram, meaning they will move monotonically in opposite senses near to the bifurcation. More precisely, suppose that there is a branch of stable equilibria $X^{s}(\lambda)$ and a branch of unstable equilibria $X^{u}(\lambda)<X^{s}(\lambda)$. Suppose that a $Y^u(\lambda)$ is the nearest branch of unstable equilibria above $X^s$. Then there is a $\lambda_1>\lambda_0$ such that
$$
\max\{ X^{u}(\mu)~:~\lambda_0<\mu<\lambda_1\} < X^{s}(\lambda)<\min\{Y^{u}(\mu)~:~\lambda_0<\mu<\lambda_1\}
$$
for all $\lambda_0<\lambda<\lambda_1$. 
The branch $X^{s}(\lambda)$ is therefore forward basin-stable and by
Theorem~\ref{thm:rtipping} case 1 there will be no R-tipping for any
parameter shifts within this range of parameters: this is a stronger conclusion than Lemma~\ref{lem:rminimal}.

\begin{figure}
\begin{center}
\centerline{\includegraphics[width=10cm]{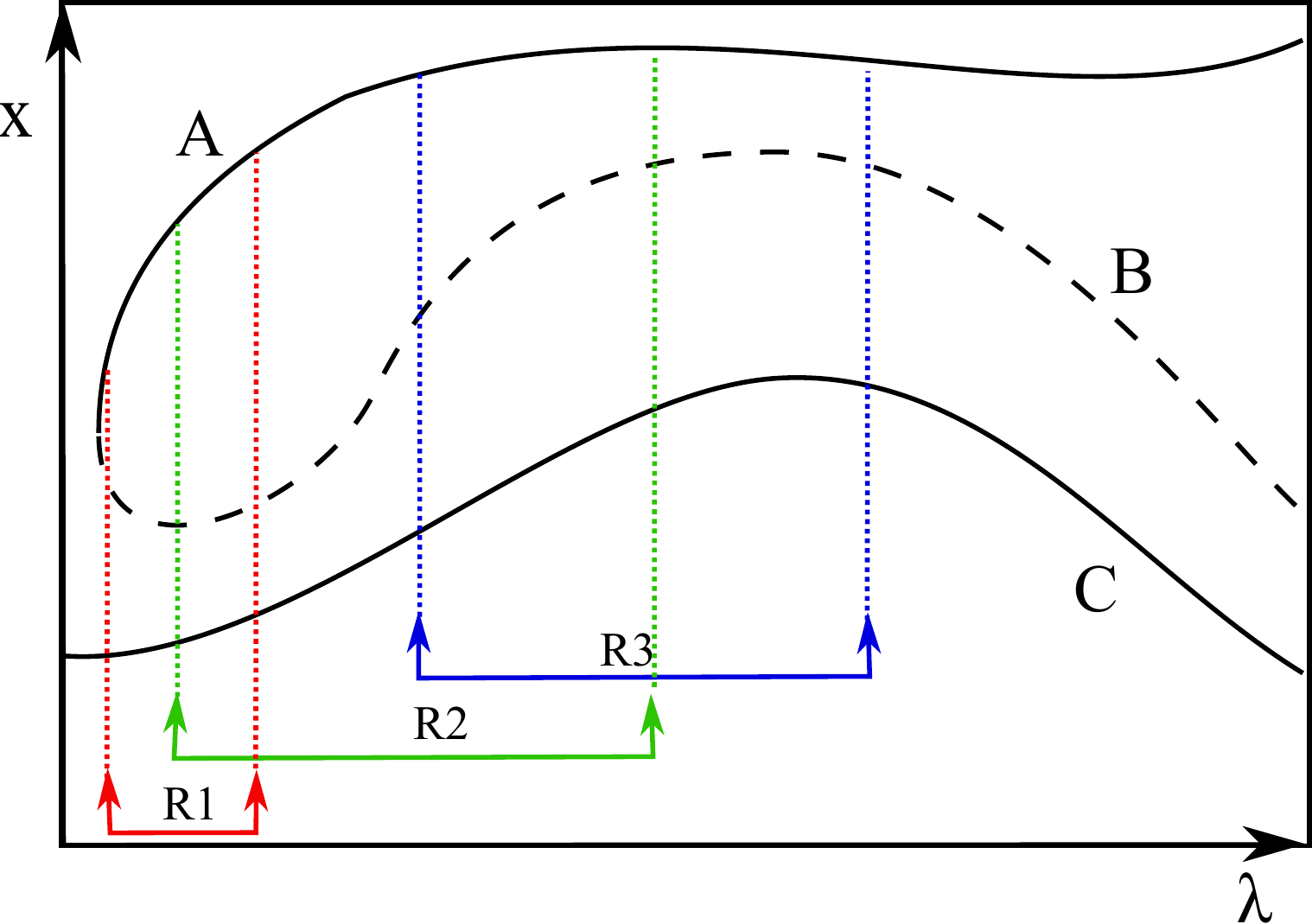}}%
\end{center}
\caption{Bifurcation diagram for a one-dimensional system with stable branches 
A and C and unstable branch B. For ranges $(\lambda_-,\lambda_+)$ with $\lambda_-<\lambda_+$ that are illustrated by R1 and R3,  Theorem~\ref{thm:rtipping} case 1 applies: the paths are forward basin stable meaning there cannot be an R-tipping from branch A. For range R2,
  Theorem~\ref{thm:rtipping} case 3 implies there can be R-tipping
  from pullback attractors starting on the stable branch A to the
  stable branch C. 
}
\label{fig:banana3}
\end{figure}

As an example see Figure~\ref{fig:banana3}; this shows a bifurcation
diagram along with three examples of choices of $(\lambda_-,\lambda_+)$. Two of these ranges (R1,R3) will not show R-tipping, while R2 can show R-tipping by Theorem~\ref{thm:rtipping} case 3.

\section{Example: R-tipping in an energy balance climate model}
\label{sec:examples}

We briefly indicate how we can characterise paths where there will be R-tipping for the Fraedrich energy-balance model of global mean temperature: see \cite{AshWieVitCox2012} for more details and an example of a path that gives $R-$tipping. We show how this can be deduced from application of Theorem~\ref{thm:rtipping}. This model for the evolution of global mean surface temperature $T$ is
\begin{equation}
  \label{eq:Fraedrich}
	\frac{dT}{dt}= f(T)=\tau^{-1}a(-T^4+b T^2-c),
\end{equation}
where the parameters $a,b,c,\tau$ can be expressed in terms of physical quantities such as the solar output, the planetary albedo etc; note that $a>0$ and $\tau>0$ from physical arguments. We set $X=T^2$ and note that (\ref{eq:Fraedrich}) transforms to
\begin{equation}
\frac{dX}{dt}= 2T\frac{dT}{dt}=\frac{2a}{\tau}(-X^2+bX-c)\sqrt{X}
\end{equation}
which (if we are bounded away from $X=T=0$) is bifurcation-equivalent to
\begin{equation}
\label{eq:be_energybal}
\frac{dX}{dt}= -X^2+bX-c
\end{equation}
in the sense of Theorem~\ref{thm:bifequiv}: though we do not need to apply this theorem for the results below. Suppose now that $b(\lambda)$ and $c(\lambda)$ vary smoothly according to some parameter $\lambda$
that undergoes a shift (without loss of generality) from $0$ to $1$. If $b(\lambda)^2>4c(\lambda)$ then there are two branches of equilibria for (\ref{eq:be_energybal}) -- and hence also for (\ref{eq:Fraedrich}); a stable branch
$$
X^{s}(\lambda)=\frac{b(\lambda)}{2}+ \frac{\sqrt{b(\lambda)^2-4c(\lambda)}}{2}
$$
and an unstable branch
$$
X^{u}(\lambda)=\frac{b(\lambda)}{2}- \frac{\sqrt{b(\lambda)^2-4c(\lambda)}}{2}.
$$
These meet at saddle nodes whenever $b(\lambda)^2=4c(\lambda)$ and there are no equilibria for $b(\lambda)^2<4c(\lambda)$. The next statement summarizes a number of results that can be proved about this system for parameter shifts of $\lambda$ from $0$ to $1$.

\begin{theorem}
\label{thm:energybal}
Consider the system (\ref{eq:Fraedrich}) with $b(\lambda)$ and $d(\lambda)$ depending smoothly on a single real parameter $\lambda\in[0,1]$. Suppose that $\lambda$ varies with time by a parameter shift $\Lambda\in \cP(0,1)$.
\begin{itemize}
\item[(a)]
If $b(\lambda)^2>4c(\lambda)$ for all $\lambda\in[0,1]$ then there is a stable branch from $X^{s}_{-}=X^{s}(0)$ to $X^{s}_{+}=X^{s}(1)$ and there is no B-tipping from $X^{s}_{-}$ for any $\Lambda\in\cP(0,1)$.
\item[(b)]
Suppose that (a) is satisfied and $X^{s}(\mu)>X^{u}(\nu)$ for all
$\mu,\nu\in[0,1]$. Then there is no irreversible R-tipping from $X^{s}_{-}$ on any $\Lambda\in\cP(0,1)$.
\item[(c)]
Suppose that (a) is satisfied and there are $0<\mu<\nu<1$  such that
$X^{s}(\mu)<X^{u}(\nu)$. Then there is a $\Lambda\in\cP(0,1)$ for which there is irreversible R-tipping from $X^{s}_{-}$.
\item[(d)]
Suppose that (a) is satisfied and let $X^{u}_{-}=X^{u}(0)$, $X^{u}_{+}=X^{u}(1)$ be the
end-points of the unstable branch. If $X^{s}_{-}<X^{u}_{+}$ then there is a
$\Lambda\in\cP(0,1)$ for which we have irreversible R-tipping from $X^{s}_{-}$.
\item[(e)]
If $b(0)^2>4c(0)$ and there is $0<\mu<1$ such that  $b(\mu)^2<4c(\mu)$ then there is B-tipping for any $\Lambda\in\cP(0,1)$ starting at $X^{s}_{-}$.
\end{itemize}
\end{theorem}

\proof Conclusion (a) follows using Lemma~\ref{lem:adiabatic}. 
Conclusion (b,c,d) follow as applications of Theorem~\ref{thm:rtipping} cases 1, 2 and 3 respectively. Conclusion (e) follows on verifying Definition~\ref{def:btipping} of B-tipping.  \qed

~

Some consequences of this are that (i) if we fix $b$ and vary just
$c$, we cannot find an R-tipping in this system: this is because in
that case
$$
X^{u}(\lambda)<\frac{b}{2}<X^{s}(\lambda)
$$
for all $\lambda\in[0,1]$ so that we can apply Theorem~\ref{thm:energybal}(b), (ii) similarly we cannot get R-tipping if we fix $c>0$ and vary just $b$. For (i) and (ii), note that it is still possible to get B-tipping in either case. (iii) Theorem~\ref{thm:energybal}(d) can be used to show that the parameter shift considered in  \cite{AshWieVitCox2012} leads to R-tipping.

\section{Discussion}
\label{sec:discuss}

In this paper we introduce a formalism that uses pullback attractors
to describe the phenomena of B-tipping and irreversible R-tipping
\cite{AshWieVitCox2012} in a class of non-autonomous systems. We
restrict to a specific type of system, namely those undergoing what we
call ``parameter shifts'' between asymptotically autonomous systems
and focus on equilibrium attractors. In this setting we show in
Theorem~\ref{thm:pullback} that there is a unique pullback point
attractor associated with each stable equilibrium for the past time
asymptotic system, and we can give definitions of what it means to
``$\epsilon$-close track'', to ``end-point track'' and to undergo some
type of ``tipping'' depending on the details of the system. For
one-dimensional systems, we obtain a number of results (in
particular Theorem~\ref{thm:rtipping}) that give necessary conditions
for irreversible R-tipping to be present or absent in a system.

We consider smooth parameter shifts that remain within the end-points
of their range and that go (possibly non-monotonically) from
$\lambda_-$ to $\lambda_+$ for some $\lambda_-<\lambda_+$. It should be
possible to relax several of these assumptions and give similar results, as long as the parameters
have well-defined limits for $t\rightarrow\pm\infty$, though this may lead to less intuitive statements.

We do not consider noise-induced tipping: for this, one needs to go beyond the
purely topological/geometric setup here, and to discuss probabilities of tipping.
We note that Theorem~\ref{thm:adiabatic} gives criteria for
tracking of stable branches by pseudo-orbits as well as useful
conditions for the tracking of branches in the presence of low
amplitude noise. Indeed, similar results can be made about the
tracking of unstable branches by pseudo-orbits though presumably if
these pseudo-orbits are the result of an occasional random process,
there will not be tracking of unstable branches for long periods of
time with any high probability. We also do not consider interaction of
noise and rate-induced tipping - but see for example \cite{Ritchie2015}.

Our results rely on an assumption of simple dynamics in state space -
in particular that the limiting behaviour of the asymptotically
autonomous systems are just equilibria, and (for the later results)
that the state space is one-dimensional and so well-ordered. 
Although many of the phenomena discussed will be present in higher
dimensional systems (for example \cite{Scheffer2008}), such systems
will clearly display additional types of nonlinear behaviour. This includes more
complex branches of attractors, periodic orbits, homoclinic and
heteroclinic connections, chaotic dynamics and complex topological
structures of attractor branches.  Nonetheless, notions such
as forward basin-stability should generalize to give similar results
for tipping from  equilibria in higher dimensions using parameter
shifts of the type considered here.

As considered in \cite{AshWieVitCox2012}, one particular type of
tipping not present in one-dimensional systems may already appear in
two dimensional state space where the system has a single globally
attracting invariant set for all values of the parameters. Such a case of transient R-tipping, associated with a novel type of
excitability threshold in
\cite{Wieczorek_etal_2010,Mitry_etal_2013,PerryWiecz2014},
cannot be explained simply in terms of the branches of the bifurcation
diagram even if we include all recurrent invariant sets.
Instead it is due to the richer geometry of attraction that can appear
in two or more dimensional systems.

We end with a brief discussion of an alternative setting that can be
used in cases where the parameter shift behaviour is determined by an
autonomous ODE, for example $d\lambda/dt=g(\lambda,r)$ for
$\lambda\in \R$ where $g(\lambda_-,r)=g(\lambda_+,r)=0$ and
$g(\lambda,r)>0$ for all $\lambda\in(\lambda_-,\lambda_+)$. In this
case we can consider an extended system of the form
$$
\frac{d}{dt}\left(\begin{array}{c} x\\ \lambda\end{array}\right) = 
\left(\begin{array}{c} f(x,\lambda) \\ g(\lambda,r) \end{array}\right).
$$
As discussed in \cite{Perryman2015}, the presence of irreversible
R-tipping can be understood in terms of (heteroclinic) connections
between saddle equilibria in this extended system and indeed this can
be used as a tool to numerically locate rates corresponding to
R-tipping. By contrast, in this paper we make no such assumption and
consider the behaviour of the non-autonomous system (\ref{eq:ode}).

\subsection*{Acknowledgements}

We thank Ulrike Feudel, Martin Rasmussen and Jan Sieber for interesting conversations and comments in relation to this work and Kate Meyer and Hassan Alkhayuon for comments that helped with the revision. We greatly thank an anonymous referee for several comments and suggesting the detailed proof for Theorem~\ref{thm:pullback} that we have included in Appendix A. We also thank Achim Ilchmann for very helpful discussions in relation to this proof. The work of PA was partially supported by EPSRC via grant EP/M008495/1 ``Research on Changes of Variability and Environmental Risk'' and partially by the European Union's Horizon 2020 research and innovation programme ``CRITICS'' under Grant Agreement number No. 643073.

\bibliographystyle{plain}

\newpage

\appendix

\section*{Appendix A: Proof of Theorem 2.2}

We thank an anonymous referee for suggesting the detailed proof that
we reproduce below.

~

\proof {\bf (of Theorem~\ref{thm:pullback})} Firstly, we fix $r>0$ and claim there is a unique solution $x(t)$ of (\ref{eq:ode}) with $x(t)\rightarrow X_-$ as $t\rightarrow -\infty$. To see this, let us define
\begin{eqnarray*}
	\omega(\epsilon)&:=& \sup \left\{|df(x,\Lambda(rt))-df(X_-,\Lambda(rt))|~:~t\in \R,|x-X_-|<\epsilon\right\}\\
	\delta(T) &:=& \max \left\{\sup_{t<-T}|f(X_-,\Lambda(rt))|,\sup_{t<-T}|df(X_-,\Lambda(rt))-df(X_-,\lambda_-)|\right\}
\end{eqnarray*}
and note, because $f\in C^2$ and $f(X_-,\lambda_-)=0$, that  $\omega(\epsilon)\rightarrow 0$ as $\epsilon\rightarrow 0$ and $\delta(T)\rightarrow 0$ as $T\rightarrow \infty$. Linear stability means that the eigenvalues of 
$$
A:=df(X_-,\lambda_-)
$$
have negative real parts and so there are $K>0$ and $\alpha>0$ such that $|e^{tA}|\leq K e^{-\alpha t}$ for $t\geq 0$ (in fact, one can choose any $\alpha$ such that $\max \{ \re(\sigma(A))\}<-\alpha<0)$: see for example \cite[Lemma~3.3.19]{HinPri2011}).

Now set $x=X_-+y$ and note that
\begin{equation}
\dot{y}= A y + h(y,t)
\label{eq:y}
\end{equation}
where
$$
h(y,t):=f (X_-+y,\Lambda(rt))-Ay
$$
so that 
\begin{align*}
dh(y,t)=&df(X_-+y,\Lambda(rt))-A\\
=&[df(X_-+y,\Lambda(rt))-df(X_-,\Lambda(rt)))]+[df(X_-,\Lambda(rt)))-df(X_-,\lambda_-)]\\
&+df(X_-,\lambda_-)-A\\
=&[df(X_-+y,\Lambda(rt))-df(X_-,\Lambda(rt)))]+[df(X_-,\Lambda(rt)))-df(X_-,\lambda_-)].
\end{align*}
Hence, for all $t<-T$ we have
$$
|h(0,t)|\leq \delta(T),~~|dh(y,t)|\leq \omega(|y|)+\delta(T)
$$
and we seek a solution of (\ref{eq:y}) such that $|y(t)|\leq \Delta$ for $t\leq -T$ for any $T>T_0$. We choose $\Delta>0$ and $T_0>0$ such that
$$
4 K\alpha^{-1}\omega(\Delta)\leq 1,~~2K\alpha^{-1}\delta(T_0)\leq \Delta,~~\mbox{ and }~4 K \alpha^{-1}\delta(T_0)\leq 1.
$$
For any continuous $y(t)$ defined for $t<-T_0$ and $|y(t)|\leq \Delta$ we define $\hat{y}$ by
$$
\hat{y}(t)=\int_{u=-\infty}^{t} e^{(t-u)A} h(y(u),u)\,du
$$
and note that the required solution of (\ref{eq:y}) satisfies $\hat{y}=y$. Note that
\begin{eqnarray}
|\hat{y}(t)|&\leq& \int_{u=-\infty}^{t} Ke^{-\alpha(t-u)} [\delta(T)+(\omega(\Delta)+\delta(T))|y(u)|]\,du\nonumber \\
& \leq& K\alpha^{-1} [\delta(T)+(\omega(\Delta)+\delta(T))\Delta].
\label{eq:yhat}
\end{eqnarray}
and hence $\|\hat{y}\|:=\sup_{t\leq -T} |y(t)|\leq \Delta$.

If there are two such functions $y_{1,2}(t)$ one can verify that
\begin{eqnarray*}
	\|\hat{y}_1-\hat{y}_2\|\leq \frac{1}{2} \|y_1-y_2\|
\end{eqnarray*}
and so the integral operator is a contraction on the set of continuous $y(t)$ with $|y(t)|\leq\Delta$ for $t<-T$ for the norm $\|\cdot\|$. Hence there is a unique fixed point such that $|y(t)|\leq \Delta$ for $t\leq -T$ and this is a solution of (\ref{eq:y}). Note from (\ref{eq:yhat}) with $\hat{y}=y$ that for $t\leq -T$
$$
|y(t)|\leq K\alpha^{-1}(\omega(\Delta)+\delta(T))\|y\|+K\alpha^{-1} \delta(T) \leq \frac{1}{2}\|y\|+K \alpha^{-1} \delta(T).
$$
Hence when $T\leq T_0$ we have $|y(t)|\leq 2K\alpha^{-1} \delta(T)$ for $t\leq-T$. Hence $y(t)\rightarrow 0$ as $t\rightarrow -\infty$. If we consider $U=B_{\Delta}(X_-)$, the
$\Delta$-ball around $X_-$, and for any $t<-T$ define
\begin{equation}
x_{pb}(t) := \bigcap_{s>0} \bigcup_{u>s} \Phi_{t,t-u}(U).
\label{eq:xpb}
\end{equation}
Then $x_{pb}(t)$ is the unique trajectory $X_-+y(t)$ for all $t<-T$ and we have $|x_{pb}(t)-X_-|<\Delta$ on $t<-T$ and $x_{pb}\rightarrow X_-$ as $t\rightarrow -\infty$. Note that (\ref{eq:xpb}) is well defined for all $t$ and is invariant under (\ref{eq:ode}) - hence it defines a unique trajectory for all $t$.

Moreover, there is a $\gamma>0$ such that for any trajectory $z(t)$ with $z(t)\in U$ for $s<t<-T$ we have 
$$
|z(t)-x_{pb}(t)|\leq K e^{-\gamma(t-s)}|z(s)-x_{pb}(s)|
$$
for $s\leq t\leq -T$. 

To see this, note that $|x_{pb}(t)-X_-|\leq \Delta= 2K \alpha^{-1}\delta(T)$ for $t\leq -T$, and hence
$$
|df(x_{pb}(t),\Lambda(rt))-A|\leq \omega(\Delta)+\delta(T).
$$
It follows from \cite[Prop 1]{Coppel} that if $Y(t)$ is the fundamental matrix solution of 
$$
\dot{y}=A(t)y=df(x_{pb}(t),\Lambda(rt))y
$$
then 
$$
|Y(t)Y^{-1}(s)|\leq K e^{-(\alpha-K(\omega(\Delta)+\delta(T)))(t-s)}
$$
for $s\leq t\leq -T$. If we choose $T$ large enough that
$$
\alpha-K(\omega(\Delta)+\delta(T))>0
$$
and put $x=x_{pb}(t)+y$ then (\ref{eq:ode}) can be written
\begin{equation}
\dot{y}=A(t)y+h_1(y,t)
\label{eq:yy}
\end{equation}
where
$$
h_1(y,t):=f (x_{pb}(t)+y,\Lambda(rt))-f(x_{pb}(t),\Lambda(rt))-df(x_{pb}(t),\Lambda(rt))y.
$$
and so if $|y|<\epsilon$ then 
$$
|h_1(y,t)|\leq \omega_1(\epsilon)|y|
$$
where
$$
\omega_1(\epsilon):= \sup\{|df(x+y,\Lambda(s))-df(x,\Lambda(s))|~:~s\in\R,~|x-X_-|\leq \Delta,|y|\leq \epsilon\}.
$$
Uniform continuity of $df$ in this region implies that $\omega_1\rightarrow 0$ as $\epsilon\rightarrow 0$. If we fix $\epsilon$ small enough that 
$$
K\omega_1(\epsilon)<\alpha-K(\omega(\Delta)+\delta(T))
$$
then it follows from the argument in \cite[Thm 9]{Coppel2} that if $y(t)$ is a solution of (\ref{eq:yy}) with $|y(s)|<\epsilon/K$ for some $s<-T$ then
$$
|y(t)|\leq K|y(s)|e^{-\gamma (t-s)}
$$
for all $t$ with $s\leq t\leq -T$ where 
$\gamma=\alpha-K(\omega(\Delta)+\delta(T))-K\omega_1(\epsilon)$. Hence $|z(t)-x_{pb}(t)|\rightarrow 0$ as $t-s\rightarrow \infty$.

Moreover, note that if $y(s)$ is a solution of (\ref{eq:yy}) where $|y(s)|<\epsilon/K$ for all $s<-T$ we have 
$$
|y(t)|\leq K|y(s)|e^{-\gamma (t-s)}
$$
for all $s,t$ with $s\leq t\leq -T$. Fixing any $t$ and taking the limit $s\rightarrow -\infty$ implies that $y(t)=0$. Hence the only solution that satisfies $|x(t)-x_{pb}(t)|<\epsilon/K$ for all $t<-T$ is the pullback attractor $x_{pb}$.
\qed

\end{document}